%% file: main.tex
\documentclass[
    english,11pt,
]{article}
\newcommand*{\Root}{.}
\newcommand*{\Dir}{.}
\long\def \forcecommand #1{\providecommand{#1}{}\renewcommand{#1}}

\newcommand*{\imageformat}{jpeg}

\usepackage{newclude}
\usepackage{\Root/src/packages}
\usepackage{\Root/src/plain-preamble}
\usepackage{\Root/src/macros}
\usepackage{\Root/src/fontconfig}

\renewcommand{\thesubfigure}{\alph{subfigure}}
\makeatletter
\renewcommand\p@subfigure{\thefigure\,({\thesubfigure})\@gobble}
\makeatother

\usepackage[normalem]{ulem}
\usepackage{authblk}

\title{\scshape \Large Repelling-Attracting Hamiltonian Monte Carlo}
\author[$^1$]{Siddharth Vishwanath%
\thanks{%
Corresponding author: %
\href{mailto:svishwanath@ucsd.edu}{\texttt{svishwanath@ucsd.edu}}%
}}
\author[$^{2,3,4}$]{Hyungsuk Tak}
\affil[$^{1}$]{Department of Mathematics, University of California, San Diego}
\affil[$^{2}$]{Department of Statistics, Pennsylvania State University}
\affil[$^{3}$]{Department of Astronomy and Astrophysics, Pennsylvania State University}
\affil[$^{4}$]{Institute for Computational and Data Sciences,  Pennsylvania State University}
\date{}

\begin{document}
\begingroup
\maketitle
\let\clearpage\relax
\include*{\Dir/inputs/abstract}%
\include*{\Dir/inputs/introduction}
\include*{\Dir/inputs/background}%
\include*{\Dir/inputs/results}%
% \include*{\Dir/inputs/experiments}%
\include*{\Dir/inputs/expt}%
\include*{\Dir/inputs/discussion}%
\include*{\Dir/inputs/proofs}%

\section*{Acknowledgement}
SV and HT would like to thank Xiao-Li Meng and David van Dyk for their insightful comments on the first draft of this manuscript. SV would also like to thank Ann Johnston for helpful discussions and feedback.

\clearpage
\bibliographystyle{abbrvnat}
\bibliography{\Dir/main}%
\endgroup
\end{document}

%% file: inputs/abstract.tex
% !TEX root = %WORKSPACE_FOLDER%/plain.tex
\begin{abstract}
    We propose a variant of Hamiltonian Monte Carlo (HMC), called the Repelling-Attracting Hamiltonian Monte Carlo (\haram{}),  for sampling from multimodal distributions. The key idea that underpins \haram{} is a departure from the conservative dynamics of Hamiltonian systems, which form the basis of traditional HMC, and turning instead to the dissipative dynamics of conformal Hamiltonian systems. In particular, \haram{} involves two stages: a mode-repelling stage to encourage the sampler to move away from regions of high probability density; and, a mode-attracting stage, which facilitates the sampler to find and settle near alternative modes. We achieve this by introducing just one additional tuning parameter---the coefficient of friction. The proposed method adapts to the geometry of the target distribution, e.g., modes and density ridges, and can generate proposals that cross low-probability barriers with little to no computational overhead in comparison to traditional HMC. Notably, \haram{} requires no additional information about the target distribution or memory of previously visited modes. We establish the theoretical basis for \haram{}, and we discuss repelling-attracting extensions to several variants of HMC in literature. Finally, we provide a tuning-free implementation via dual-averaging, and we demonstrate its effectiveness in sampling from, both, multimodal and unimodal distributions in high dimensions.
\end{abstract}

%% file: inputs/introduction.tex
% !TEX root = %WORKSPACE_FOLDER%/plain.tex

\section{Introduction}
\label{sec:intro}

Consider the problem of drawing samples $\qty{\q_1, \q_2, \q_3, \dots} \subset \R^d$ from a target distribution, $\pi$, whose probability density function is given~by
\eq{
    \pi(\q) = \f 1Z \exp( - U(\q) ),\nn
}
where $\q=({q_1, q_2, \ldots, q_d}) \in \R^d$, $U(\q) := -\log\pi(\q)$, and $Z = \int \exp(-U(\q))d\q$ is its normalizing constant. Such problems are prevalent in Bayesian inference (\citealp{gelman1995bayesian}), probabilistic machine learning (\citealp{andrieu2003introduction}), and statistical physics (\citealp{binder1992monte}), among other fields, and constitutes the core of Markov chain Monte Carlo (MCMC) methods. The goal of MCMC is to generate a sequence of dependent samples $\qty{\q_1, \q_2, \dots, \q_n}$ that converges to the target distribution $\pi$ as $n \to \infty$ \citep{brooks2011handbook}.
% is the log-pdf of $\pi$

In most practical settings, the analytical expression for the normalizing constant $Z$ is seldom available, or, its numerical computation is intractable. Therefore, sampling from the target distribution $\pi(\q)$ can become an insurmountable task \citep{moller2006efficient,murray2012mcmc,park2018bayes}. The Metropolis-Hastings algorithm \citep{metropolis1953equation, hastings1970monte}, overcomes this by generating proposals from a random walk. However, the random walk nature of the Metropolis-Hastings algorithm makes it inefficient for sampling from high-dimensional target distributions.

\begin{figure}[t!]
    \centering
     \includegraphics[width=0.6\textwidth]{\Root/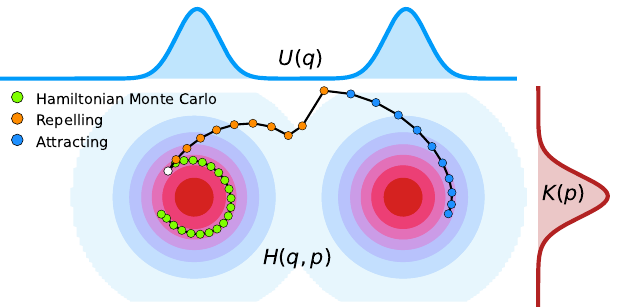}
    \caption{Illustration of \haram{} and HMC transitions for a bimodal target distribution $\pi(q)=-\log(U(q))$ in an extended phase space; $U(q)$ is the potential energy, $K(p)$ is the kinetic energy, and  $H(q, p)=U(q)+K(p)$ is the   Hamiltonian. The \haram{} trajectory concatenates two conformal Hamiltonian trajectories to generate the final proposal. The first trajectory (in orange) simulates the mode-repelling dynamics, forcing the trajectory to move away from the current mode. The second trajectory (in blue) simulates the mode-attracting dynamics, and the trajectory is drawn towards a different mode. Combined, the \haram{} trajectory is able to jump between modes by moving \textit{across} the energy level sets of $H(q, p)$. This is in contrast to HMC, which moves \textit{along} the energy level sets of $H(q, p)$, and is unable to jump between modes.}
    \label{fig:idea}
\end{figure}

Hamiltonian Monte Carlo (HMC; \citealp{duane1987hybrid,neal2011mcmc}) along with its tuning-free implementations, such as the No-U-Turn sampler (NUTS; \citealp{hoffman14a}), have emerged as the mainstay for sampling from high dimensional distributions when the gradient information, $\grad\log\pi(\q)$, is available. HMC broadly falls under the category of MCMC methods based on data augmentation \citep{tanner1987calculation, dyk2001art}, whereby an auxiliary variable $\p=({p_1, p_2, \ldots, p_d}) \in \R^d$, specified by its density function $\pi(\p) \propto \exp(-K(\p))$, is introduced such that the joint distribution of $(\qp)$~is~given~by
\eq{
    \pi(\qp) \propto \exp\qty\big( - U(\q) - K(\p)  ) =: \exp( -H(\qp) ).\nn
}
By regarding $\q$ as its position of a particle with potential energy $U(\q)$, and $\p$ as its conjugate momentum with kinetic energy $K(\p)$, the term $H(\q, \p)$ can be interpreted as the Hamiltonian of a physical system in an extended phase space $(\qp) \in \R^{2d}$. Given the current state $\q$, the key idea of HMC is to sample ${\p \sim \pi(\p)}$ from a known distribution and generate a proposal $(\qp) \mapsto (\qp[T])$ by simulating the Hamiltonian dynamics of the system for a fixed time $T$, i.e.,
$$
\dd{t}\q_t = \D_{\p} H(\p_t, \p_t) = \grad K(\p_t), \qq{and} \dd{t}\p_t = -\D_\q H(\q, \p) = \grad U(\q).
$$
The green trajectory in Figure~\ref{fig:idea} illustrates a typical HMC proposal. Unlike the random-walk Metropolis method, by leveraging information from $\grad U(\q) = -\grad\log\pi(\q)$, HMC is able to adapt to the geometry of the target distribution, and produce distant proposals with high acceptance probabilities even in high dimensions. We refer the reader to \cite{neal2011mcmc} and \cite{betancourt2017conceptual} for a detailed review of HMC.

%for the same number of leapfrog steps (that is, with comparable computational cost).

Several lines of work have extended the HMC algorithm to situations where the $\log\pi(\q)$ is discontinuous \citep{nishimura2020discontinuous}, when the support of $\pi$ is discrete \citep{zhang2012continuous}, when the support of $\pi$ is a Riemannian manifold \citep{girolami2011rhmc}, or when the gradient evaluation is intractible \citep{strathmann2015gradient}. However, a fundamental shortcoming of HMC is its inability to sample from multimodal distributions, and is widely acknowledged in literature (e.g., \citealp{barp2018geometry}). In particular, since $\grad U(\q)$ only captures local information, HMC is unable to traverse and generate proposals across well-seaprated regions of $\pi$. This is particularly problematic in high-dimensional settings, where zero-density regions are more prevalent, and the modes of $\pi$ are more likely to be well-separated.

In order to adress the shortcomings of HMC for multimodal target distributions, several methods have been proposed in the literature. Broadly, they may be classified into the following three categories: (i)~methods that generalize the original HMC dynamics to that of another physical system \citep[e.g.,][]{leimkuhler2009metropolis,tripuraneni2017magnetic,Liu2019QuantumInspiredHM}, (ii)~methods which lower the energy barrier between modes by augmenting the phase-space via auxiliary variables, or via tempering schemes \citep[e.g.,][]{betancourt2014adiabatic,nishimura2016geometrically,graham2017continuously,nemeth2019pseudo,park2021sampling,alenlov2021pseudo}, and (iii)~methods that introduce a systematic method for identifying where the modes of the distribution are, and then adapting the HMC dynamics to jump between these modes \citep[e.g.,][]{sminchisescu2011generalized,wormhole2014,pompe2020framework}. Section~\ref{sec:multimodal} provides a more detailed review of these methods. Implementing these methods often requires additional tuning parameters, and are computationally more expensive than traditional HMC. 

In this paper we propose the Repelling-Attracting Hamiltonian Monte Carlo (\haram{}) algorithm for sampling from both unimodal and multimodal distributions. The philosophy underlying \haram{} is to generate proposals in two stages: in the first stage, high-density regions of $\pi$ act as repelling states and encourage the trajectory to move away from them into low density regions, and, in the second state, the high-density regions act as attracting states and encourage the trajectory to settle near them. The orange and blue trajectories in Figure~\ref{fig:idea} illustrate a typical \haram{} proposal. This is achieved by moving away from the physical interpretation of Hamiltonian dynamics and by considering, instead, conformal Hamiltonian systems where energy can be created and dissipated. Therefore, our method can be described as a middle-ground between the perspectives, (ii) and (iii), described above. Our main contributions are the following:

\begin{itemize}[ leftmargin=!, itemindent=1.5em, itemsep=0em]
\item We propose a new MCMC sampler, \haram{}, that is designed to be a versatile sampler for both unimodal and multimodal distributions, which requires only one additional tuning parameter, $\gamma$, the coefficient of friction, in addition to the default tuning parameters of HMC.
\item We show that \haram{} satisfies the key properties which make HMC a successful MCMC sampler, and demonstrate that \haram{} can maintain high acceptance rate across long trajectories across multiple modes even in high dimensions.
\item We provide an automatic tuning procedure of \haram{} using dual-averaging method of \cite{nesterov2009dual}, akin to the automatic tuning procedure for the No-U-Turn sampler and HMC \citep{hoffman14a}.
\item We show that \haram{} can be generalized to incorporate other variants of HMC, e.g., relativistic Monte-Carlo \citep{lu2017relativistic}, magnetic HMC \citep{tripuraneni2017magnetic}, and the non-canonical HMC \citep{brofos2020non}, among others. 
\end{itemize}

The code for the implementation and for reproducing the experiments is made publicly available at \linebreak{\footnotesize {\url{https://github.com/sidv23/ra-hmc}}}.

\bigskip

\textbf{Notations.}\quad For a space $X$, the map $\id: X \to X$ is the identity map given by $\id(x) = x$ for all $x \in X$. We use $\onev_d, \zerov_d \in \R^d$ to denote the vector of ones and zeros in $\R^d$, $\I_d$ to denote the $d \times d$ identity matrix, and $\O_d$ to denote the $d \times d$ zero matrix. We drop the subscript $d$ when the dimensionality is clear from the context. For a vector $\xv \in \R^{2k}$, we use $\textup{proj}_1(\xv)$ and $\textup{proj}_2(\xv)$ to be the projection of $\xv$ onto the first and last $k$ components, respectively. $\calo{N}(\mu, \S)$ denotes the multivariate Gaussian distribution with mean $\mu \in \R^d$ and covariance matrix $\S \in \R^{d \times d}$. For a differentiable function $f:\R^d \to \R$ the gradient and Hessian of $f$ are denoted by $\grad f$ and $\grad^2 f$, respectively, and for $F: \R^d \rightarrow \R^p$, $\jac{F}$ is  the $d \times p$~Jacobian~matrix.\newline As per convention in MCMC literature, and with a slight overload of notation, we use $\pi$ as a general notation to represent the target distribution that we aim to sample from, and $\supp(\pi)$ to denote its support. The exact distribution will be made precise in context. Throughout this work, we assume that --$\log\pi(\q)$ is continously differentiable and $\ell$-Lipschitz for some $\ell > 0$. 

%% file: inputs/background.tex
\section{Background}
\label{sec:background}
In this section, we provide the necessary background on Hamiltonian Monte Carlo (HMC), which forms the basis of our proposed method. We begin by detailing the general methodology of HMC in Section~\ref{sec:hmc}, and summarize the key characteristics that contribute to the success of HMC in Section~\ref{sec:properties}. Finally, in Section~\ref{sec:multimodal}, we explore the primary challenges when sampling from multimodal distributions and examine the existing literature that tackles this issue. We refer the reader to \cite{neal2011mcmc,betancourt2017conceptual} and \cite{betancourt2017geometric} for a comprehensive overview of HMC and its theoretical underpinnings.

\subsection{Hamiltonian Monte Carlo}
\label{sec:hmc}

Let $\pi(\q)$ be a target distribution whose probability density function is given by $\pi(\q) \propto \exp(-U(\q))$, where $\q=({q_1, q_2, \ldots, q_d}) \in \R^d$. By augmenting the target distribution with an auxiliary momentum variable $\p=({p_1, p_2, \ldots, p_d}) \in \R^d$, whose marginal distribution is given by $\pi(\p) \propto \exp(-K(\p))$ for a user-specified choice of $K(\p)$, the joint target distribution of $(\qp)$ is given by
\eq{
    \pi(\qp) \propto \exp\qty\big( - U(\q) - K(\p)  ) =: \exp( -H(\qp) ).\nn
}
With an appropriate distribution for $\p$, the resulting energy function ${H(\qp) := U(\q) + K(\p)}$ can be characterized as the Hamiltonian for the system in an \mbox{\textit{extended phase space}} $(\qp) \in \R^{2d}$, where $\q$ corresponds to the position of a particle, and $\p$ corresponds to its (conjugate) momentum. The functions $U(\q) := -\log\pi(\q)$ and $K(\p)$ are commonly referred to as the potential energy function and the kinetic energy functions, respectively. 

Hamiltonian Monte Carlo (HMC) generates samples from $\pi(\q)$ by first resampling auxiliary momentum $\p \sim \pi(\p)$, and then simulating Hamiltonian dynamics in the extended phase space $(\qp) \in \R^{2d}$ for a time interval $T$. The auxiliary momentum $\p$ is subsequently discarded to obtain samples from the target distribution $\pi(\q)$. The HMC algorithm proceeds as follows: given the initial Markov chain state $(\qp[0])$, the particle is subject to the dynamics given by
\eq{
    \dd{t} \vec{\q_t, \p_t} = \underbrace{\mat{\O & \I \\-\I & \O}}_{=:\Om} \vec{ \grad U(\q_t), \grad K(\p_t)  }.
    \label{eq:hamflow-long}
}
\vspace*{-0.5em}
By taking $\z_t = \pa{\qp[t]}$, \eref{eq:hamflow-long} can be concisely written as the Hamiltonian flow
\eq{
    \dd{t} \z_t = \Om \grad_{\!\z} H(\z_t),
    \label{eq:hamflow}
}
where the skew-symmetric matrix $\Om$ can be thought of as the coordinate representation for the symplectic structure underlying the Hamiltonian system. The position of the particle at time~$t$ is given by the integral operator $\P_t: \Rdd \rightarrow \Rdd$, where $\z_{t} = \P_t(\z_0)$ is obtained by solving the autonomous system of differential equations in \eref{eq:hamflow}. Then, the proposed state at time $T$ is obtained by applying a momentum-flip operation to $\z_T$, that is, $\tz_T = \flip(\z_T)$ where ${\flip: (\qp) \mapsto (\q, -\p)}$. Finally, the proposed state $\tz_T$ is accepted/rejected via a Metropolis-Hastings acceptance probability given the current state $\z$, 
\eq{
    \alpha( \tz_T \mid \z) := \min\qty\Bigg{1,~ \exp\qty\Big( -H(\tz_T) + H(\z) ) \times \abs{\det\pa{ \f{\d \tz_T}{\d \z} }}}.
    \label{eq:mh-ratio}
}

\begin{remark}
    The HMC methodology offers flexibility in the choice of the distribution for the auxiliary momentum $\p$, $\pi(\p) \propto \exp(-K(\p))$. However, in practice, $K(\p)$ is typically taken to~be 
    $$
        K(\p) = \f 12 \p^\top\S\inv\p,
    $$
    such that the momentum variable is sampled from a Gaussian distribution, i.e., ${\p \sim \N(\zerov, \S)}$.
\end{remark}

While the flow operator $\P_T$, corresponding to the Hamiltonian flow in \eref{eq:hamflow}, is seldom available analytically, efficient numerical integrators based on splitting methods are widely available, such as the Störmer-Verlet/leapfrog scheme \citep{hairer2006geometric}. In the leapfrog scheme, the Hamiltonian of the system is split into three parts, i.e., 
\eq{
H(\z) = H_1(\z) + H_2(\z) + H_3(\z) = \f 12 U(\q) + K(\p) + \f 12 K(\q).\nn
} 
Then, for a step-size $\e \approx dt$ and $i=1,2,3$, the flow $\P_{dt}$ associated with each sub-Hamiltonian $H_i$ in time interval $(t, t+\e]$ is numerically approximated by maps $\xi_{\e, i}$. The resulting flow is given by ${\P_\e := \xi_{1} \circ \xi_{2} \circ \xi_{3} \approx \P_{dt}}$, where
\eq{
    \xi_{1}(\qp) = \xi_{3}(\qp) =\qty\Big( \q,  \ \ \p - \f \e2 \grad U(\q) ) \quad \text{ and } \quad
    \xi_{2}(\qp) = \qty\Big( \q + \e\S\inv \p, \ \ \p ).
    \label{eq:xis}
}
Thus, the position of a particle  at time $T$, $\z_T$, is numerically approximated  by applying  $\P_\e$ consecutively $L= \floor{T/\e}$ times, i.e., 
$$
\z_T = \P_{\e, L}(\z) := \qty(\P_\e)^{L} (\z) = \P_\e \circ \P_\e \circ \cdots \circ \P_\e(\z).
$$
% obtained 

\subsection{Key Properties of Hamiltonian Monte Carlo}
\label{sec:properties}
The success of HMC as a sampling method crucially hinges on the following properties.
%is applied $L$ times

% \begin{description}[style=unboxed, leftmargin=0.1cm, itemsep=1em]
\begin{enumerate}[style=unboxed, leftmargin=0cm, itemsep=1em, label=\bfseries{$\calo{P}$\arabic*.}, ref = \bfseries{($\calo{P}$\arabic*)}]
    \item \label{p:reversible} \textbf{Reversibility.} Given the current state $\z$, the proposed states, $\z \mapsto \P_T(\z)$ and $\z \mapsto \P_{\e, L}(\z)$, in HMC constitute a deterministic Markov kernel. From \cite{tierney1994markov}, a deterministic Markov kernel ${\kappa: \Rdd \rightarrow \Rdd}$ satisfies detailed-balance if and only if $\kappa$ is an involution, i.e., $\kappa\,\circ\,\kappa = \id$. Although $\P_T$ and $\P_\el$, by themselves, are not involutions; when they are composed with the momentum-flip operator, $\flip$, the resulting transition kernels, $\flip\circ\P_t$ and $\flip \circ \P_\el$, are involutions, i.e., ${(\flip \circ \P_t) \circ (\flip \circ \P_t) = \id}$, and  ${(\flip \circ \P_{\e, L}) \circ (\flip \circ \P_{\e, L}) = \id}$. Therefore, the momentum flip operation in \eref{eq:mh-ratio} is crucial for ensuring that the resulting Markov chain satisfies the detailed-balance condition.

    \item \label{p:volume} \textbf{Volume Preservation.} Liouville's theorem \citep[Chapter~3.16~B]{arnol2013mathematical} guarantees that the Hamiltonian flow, $\P_t$, preserves volume in the phase space, i.e., $|\det(\jac{\P_t(\z)})| = 1$ for all $t \in \R$. Since the momentum flip operator $\flip$ is also volume-preserving, the Jacobian term appearing in \eref{eq:mh-ratio} becomes
    \eq{
        \det\pa{ \f{d\tz_T}{d\z} } = \det\pa{ \jac{(\flip \circ \P_T)}(\z) } = \det\pa{\flip} \times \det\pa{ \jac{\P_T}(\z) } = 1,\nn
    }
    where the second equality follows from the Jacobian chain rule. This ensures that the Jacobi determinant appearing in \eref{eq:mh-ratio} is always $1$, and substantially reduces the computational cost of computing $\alpha$ in \eref{eq:mh-ratio}. See, for example, \cite{levy2018generalizing}, where sacrificing volume preservation warrants careful adjustments in order to make the resulting algorithm computationally tractable.

    \item \label{p:energy} \textbf{Energy conservation.} Hamiltonian dynamics traverses iso-energy contours in the extended phase space, and the flow $\P_t$ preserves the Hamiltonian, i.e.,
    \eq{
        \dd{t} H(\z_t) = 0.\nn
    }
    Thus, $H\qty(\P_t(\z)) = H(\z)$ for all $t \in [0,~ T]$. Combined with the fact that the momentum-flip operator preserves the Hamiltonian, i.e., $H(\flip(\z)) = H(\z)$, it follows that  $ H(\tz_T)={H\qty( \flip \circ \P_T(\z )) = H(\z)}$, and the energy in the proposed state $\tz_T$ is also preserved. Therefore, in combination with \ref{p:volume}, if the Hamiltonian dynamics can be exactly simulated, then the resulting proposal $\tz_T$ can be accepted with probability~$1$.

    \item \label{p:symplectic} \textbf{Symplecticity.} Hamiltonian systems encode strong symmetries, which are made precise using machinery from differential geometry; see \cite{betancourt2017geometric} for an overview. In particular, the \textit{symplectic} $2$-form, $\omega$, encapsulates the conservation laws and dynamical evolution prescribed by Hamiltonian dynamics. Informally, the {symplectic} $2$-form, denoted by the wedge product 
    $$
    \omega = \textrm{d}\q \wedge \textrm{d}\p = \sum_{i=1}^d \textrm{d}q_i \wedge \textrm{d}p_i, 
    $$
    measures the oriented area spanned by infinitesimal displacements in position, $\q$, and the momentum, $\p$. Remarkably, in addition to the volume preservation in \ref{p:volume}, the Hamiltonian flow $\P_t$ preserves symplectic form $\omega$, or, equivalently, from \citet[Chapter~9.44~D]{arnol2013mathematical},
    $$
    \qty({\f{d\z_t}{d\z}})^\top \Om\inv \qty({\f{d\z_t}{d\z}}) =  \jac{\P_t}(\z)^\top \Om\inv \jac{\P_t}(\z) = \Om\inv \qq{for all $t \in \R$.}  
    $$
    Therefore, it follows that preserving symplectic structure also guarantees that volume is conserved.

    Notably, however, preserving symplectic structure provides access to \textit{symplectic integrators}, $\P_{\e, L}$, to approximate the continuous-time flow $\P_T$. Like their continuous time analogues, symplectic integrators also preserve the symplectic form (and, therefore, volume) in phase space, and satisfy reversibility \`{a}-la \ref{p:reversible} \citep{hairer2006geometric}. Some notable examples include  (first-order) symplectic Euler integrator, the (second-order) leapfrog integrators and (higher-order) Yoshida integrators. 
    
    Moreover, even though symplectic integrators do not preserve the Hamiltonian $H(\z)$, they exactly preserve a \textit{shadow Hamiltonian} \citep[Chapter IX,~Theorem~3.2]{hairer2006geometric} whose approximation error is strongly controlled for the integration time $T$. For instance, for the leapfrog integrator the global error is
    \eq{
        \norm\big{ \P_T(\z) - \P_{\e, L}(\z)} = O(\e^2) \qq{for} T = \epsilon L.\nn
    }
    Moreover, it can be shown (see, e.g., \citealp[Chapter IX,~Section~8]{hairer2006geometric}) that there exists a constant $C > 0$ such that for exponentially long simulation lengths, $T = O(\e^2e^{C/\e})$, the energy drift is controlled by $\abs{H(\P_T(\z)) - H(\P_\el(\z))} = O(\e^2)$. Consequently, the acceptance ratio for the proposed state $\tz_T=\flip \circ \P_{\e, L}(\z)$ in \eref{eq:mh-ratio} is
    \eq{
        \alpha( \flip \circ \P_{\e, L}(\z) \mid \z ) = \exp(-{O(\e^2)}), \nn
    }
    which ensures that the acceptance probability remains relatively high, even for fairly distant trajectories. 
% \end{description}
\end{enumerate}

%%%%%%%%%%%%%%%%%%%%%%%%%%%%%%%%%%%%%%%%%%%%%%%%%%%%%%%%%%%%%%%%%%%%%%%%%%%%%%%%%%%%
%%%%%% SECTION: MULTIMODAL

\subsection{Drawbacks of Hamiltonian Monte Carlo for Multimodal Target Distributions}
\label{sec:multimodal}

The inability of HMC to sample effectively from multimodal distribution has been widely acknowledged in literature (e.g., \citealp{celeux2000computational, sminchisescu2011generalized}). For example, \cite{mangoubi2018does} show that, for multimodal target distributions, the ability of HMC to transition between modes is provably worse than that of random-walk Metropolis. We begin by providing an illustrative example that demonstrates the difficulty of HMC in exploring a multimodal distribution when $\pi$ is a mixture of two Gaussian distributons.

\begin{proposition}
    For $d \ge 2$, let $\pi \propto 0.5\,\N(-\bd, \I_d) + 0.5\,\N(\bd, \I_d)$ be a bimodal target distribution, where $b > 0$. Let $\eta(t; \zo) = \textup{proj}_1 \circ \P_t(\zo)$, i.e., $\eta(t; \qo, \po) = \q_t$ where $(\qt, \pt) = \P_t(\qo, \po)$, and let $\Aa$ be the the high-density region around the the mode $\bd$ given by
    $$
        \Aa = \qty{ \q \in \R^d : \norm{\q - \bd} \le \alpha \sqrt{d}}.
    $$
    Define $\E(b, d, \alpha)$ to be the event that the HMC transition kernel proposes a point $\eta(T; \qo, \po)$ which is closer to $-\bd$ than $\bd$, i.e.,  
    \eq{
        \E(b, d, \alpha) = \qty\Big{ \norm{\eta(T; \qo, \po) + \bd} \le \norm{\eta(T; \qo, \po) - \bd}, \qo \in \Aa }.\nn
    }    
    Then, for $b \ge \sqrt{\alpha^2 + 2 \sigma^2}$,
    \eq{
        \pr\qty\Big( \E(b, d, \alpha) ) \le \exp\qty( -\frac{d}{2} \qty{ \frac{b^2 - \alpha^2}{2\sigma^2} - 1 - \log\qty(\frac{b^2 - \alpha^2}{2\sigma^2})} ).
        \label{eq:mode-transition}
    }
    \label{prop:mode-transition}
\end{proposition}

The proof of Proposition~\ref{prop:mode-transition} is provided in Appendix~\ref{proof:prop:mode-transition}, and we make the following observations from the result.

\begin{remark}
    Although just an upper bound, \eref{eq:mode-transition} is useful for understanding the exponential decay of the transition probability as a function of the separation $b$ and the dimension $d$.

    \begin{enumerate}[label=\textup{(\roman*)}, align=left, leftmargin=!, itemindent=0.5em]
        \item For simplicity, if we take $\alpha^2 = (1-\delta)b^2$ for sufficiently small $\delta \in (0, 1)$, then the resulting term, $\delta b^2 / \sigma^2$ in \eref{eq:mode-transition}, effectively quantifies the energy barrier HMC needs to overcome to transition from the high-density region $\Aa$ near mode $\bd$, and propose a state closer to the mode $-\bd$. By applying the inequality, $x - 1 - \log(x) \ge (x-1)^2/x$, it follows that the mode-transition probability bound in \eref{eq:mode-transition} decays exponentially with an increase in the energy barrier and the dimension $d$.
        \item The condition that $b \ge \sqrt{\alpha^2 + 2\sigma^2}$ requires that the modes of $\pi$ are sufficiently well-separated for the bound in Proposition~\ref{prop:mode-transition} to hold. Indeed, if $b$ was small, then the effective energy barriers between the modes would be small, and the effect of the multimodality of the target distribution would be negligible. 
        \item We note that it may be possibe to obtain bounds similar to \eref{eq:mode-transition} for more general multimodal distributions by studying the conductance of the HMC algorithm for the target distribution $\pi$ \citep{mangoubi2018does}. We also point out that the bound in \eref{eq:mode-transition} is different from Eq.~(5) of \cite{park2021sampling}, which only provides a bound for transitioning between modes for a particular path in the HMC trajectory starting at one mode of the target distribution $\pi$. 
    \end{enumerate}
\end{remark}

In order to adress these shortcomings of HMC for multimodal target distributions, several methods have been proposed in the literature. As noted in Section~\ref{sec:intro}, they can be, broadly, classified into the following three categories: (i)~methods that generalize the original HMC dynamics to that of another physical system, (ii)~methods which lower the energy barrier between modes by augmenting the phase-space via tempering schemes or by introducing auxiliary variables, and (iii)~methods that introduce a systematic method for identifying where the modes of the distribution are, and then adapting the HMC dynamics to jump between these modes.

In the first category, \cite{tripuraneni2017magnetic} propose magnetic HMC where the HMC dynamics from classical mechanics is replaced by the Newtonian dynamics of a charged particle in a magnetic field. In essence, their method considers a non-canonical symplectic structure by replacing the lower-diagonal block of $\Om$ in \eref{eq:hamflow} with a user-specified invertible matrix $\mathbf{R}$. \cite{brofos2020non} extend this method to general non-canonical Hamiltonian systems by transforming the original phase space using the Darboux basis \citep[Chapter~8.43~B]{arnol2013mathematical}. Although these methods ensure that the resulting dynamics preserve \ref{p:reversible}--\ref{p:symplectic},  an appropriate  choice of the matrices for the non-canonical Hamiltonian system remains an open question, and may be largely problem-specific.

In the second category, \cite{betancourt2014adiabatic} proposes adiabatic Monte Carlo based on the principles of thermodynamics by replacing the symplectic structure with a contact structure \cite[Appendix~4]{arnol2013mathematical}. This results in a \textit{contact Hamiltonian} defined on a $(2d+1)$-dimensional phase space given by
\eq{
    H_C(\q, \p, \eta) = \beta(\eta) U(\q) + K(\p) + (1 - \beta(\eta))V(\q) + \log Z\qty(\beta\qty(\eta)), \nn
}
where $\eta\in \R$ governs a continuously varying (inverse) temperature $\beta(\eta) \in [0,1]$. Roughly speaking, this temperature interpolates between the target $\pi$ and a user-specified reference distribution $\varpi(\q) \propto \exp(-V(\q))$. Although theoretically appealing, practical implementation of the adiabatic Monte Carlo may be somewhat esoteric, and turns out to be infeasible in cases where the log-partition function $Z(\beta) = \int\pi(\q)^\beta d\q$ is not known. \cite{graham2017continuously} and \cite{luo2018thermostat} show that this issue can be, to some extent, ameliorated by using the approach of \cite{gobbo2015extended}, i.e., by lifting the contact Hamiltonian to an \textit{extended Hamiltonian} in a $(2d + 2)-$dimensional phase space. Notwithstanding, an appropriate choice of the confining potential, $V$, and a suitable choice of the ``conjugate momentum'' for $\beta$ are not entirely obvious. \cite{park2021sampling}, on the other hand, proposes tempered HMC based on a pre-specified mass scaling schedule $\qty{\S_t : t \in \R}$ for the mass matrix along the flow $\P_t$. While the method performs well in exploring multimodal distributions, the time-varying mass schedule requires additional tuning, and, notably, doesn't maintain the symplectic structure. On a different note, \cite{nemeth2019pseudo} propose the pseudo-extended HMC where the target distribution is a product of the original $\pi(\q)$ and its tempered versions $\qty{\pi(\q)^{\beta_1}, \dots, \pi(\q)^{\beta_k}}$, with the expectation that in the augmented phase-space $\R^{2dk}$, the modes are connected to each other and facilitates better mode exploration for HMC. 

The third class of HMC variants for sampling from multimodal distributions carefully design methods for identifying the location of modes while sampling at the same time. The difference between these methods, broadly, lies in the way they identify the modes (inner optimization routines) and how this information is transmitted to the Markov chain which is sampling from the target distribution. For example, \cite{sminchisescu2011generalized} propose the  method in which the task of finding the modes and sampling within modes is integrated into a single framework. \cite{wormhole2014} propose the Wormhole Monte Carlo whereby, instead of performing local moves, mode-transition is enabled by constructing a special Riemannian metric which accounts for the geometry of the original phase space \textit{and} the location of the modes. Their method uses the Riemann manifold Hamiltonian Monte Carlo  algorithm \citep{girolami2011rhmc} to generate proposals. %\cite{pompe2020framework} propose the jumping adaptive multimodal sampler in which the task of finding the modes and sampling within modes is integrated into a single framework while dynamically learning the optimal parameters during runtime. \cite{tawn2021annealed} propose the annealed leap-point sampler where annealing is used, instead of tempering, to explore the modes of the target distribution.

Although these methods excel at investigating multimodal distributions due to their specific design to overcome the limitations of Hamiltonian Monte Carlo (HMC), they come with a trade-off. Often, the additional tuning parameters and  computational overhead they  introduce might be unnecessary and harder to justify when sampling from simpler, unimodal distributions. As we point out in Section~\ref{exp:unimodal}, our proposed methodology demonstrates advantages for sampling from both unimodal and multimodal distributions, and and requires just one additional tuning parameter in comparison to HMC.

%% file: inputs/results.tex
% !TEX root = %WORKSPACE_FOLDER%/plain.tex

%%%%%%%%%%%%%%%%%%%%%%%%%%%%%%%%%%%%%%%C%%%%%%%%%%%%%%%%%%%%%%%%%%%%%%%%%%%%%%%%%%%%%
%%%%%% SECTION: Friction

\section{Repelling-Attracting  Hamiltonian Monte Carlo}\label{sec:raHMC}

With this background, in this section we introduce the proposed repelling-attracting Hamiltonian Monte Carlo (\haram{}) algorithm. We first describe the key ideas behind \haram{} in Section~\ref{sec:motivation}, and, in Section~\ref{sec:raHMCdetail}, we describe its implementation details and investigate its properties as a Markov chain Monte Carlo algorithm. Next, in Section~\ref{sec:tuning} we explain the procedure for tuning the parameters of \haram{}, and briefly discuss extensions of \haram{} in Section~\ref{sec:extension}.

\subsection{Motivation}
\label{sec:motivation}

The key idea behind \haram{} stems from non-idealized Hamiltonian systems in the presence of dissipative, damping forces, e.g., friction, viscosity, etc., termed as \textit{conformal Hamiltonian systems} \citep{mclachlan2001conformal}.

%%%%%%%%%%%%%%%%%%%%%%%
% \subsection{Friction in Hamiltonian dynamics}
%%%%%%%%%%%%%%%%%%%%%%%

{\bfseries Adding friction dissipates energy.\quad} The dynamics of a particle in a linearly damped  conformal Hamiltonian system is governed by the equation
\eq{
    \dd{t} \vec{\q_t, \p_t} = \mat{\O & \I \\-\I & \O} \vec{ \grad U(\q_t), \grad K(\p_t)  } - \mat{\O & \O \\ \O & \gamma\I} \vec{ \q_t, \p_t },
    \label{eq:conf-hamflow-long}
}
where $\gamma \ge 0$ is a damping constant, which, in the case of classical mechanical systems, is also called the \textit{coefficient of friction}. Similar to \eref{eq:hamflow}, by taking $\z_t = (\qp[t])$ to represent the state of the particle in the extended phase space, we can succinctly rewrite \eref{eq:conf-hamflow-long} as
\eq{
    \dd{t}\z_t = \Om \grad_{\!\z} H(\z_t) - \G \z_t, \qq{where} \G:=\mat{\O & \O \\ \O & \gamma\I}.
    \label{eq:conf-hamflow}
}
When $\gamma=0$, the dynamics of the conformal Hamiltonian system reduces to that of the conservative Hamiltonian system in \eref{eq:hamflow-long}. On the other hand, when $\gamma>0$, the friction inhibits the motion of a particle since $\d\p_t/\d t= \grad U(\q_t) - \gamma \p_t$. In contrast to the interpretation of conservative Hamiltonian dynamics in \citet{neal1996hmc} as that of a hockey puck sliding over a frictionless surface, the conformal Hamiltonian dynamics in \eref{eq:conf-hamflow} can be interpreted as a dynamics associated with the same puck on a rough surface, thereby dissipating energy. 

The dynamics of a particle under a conformal Hamiltonian vector field encourages the particle to \textit{flow towards} the critical points of $U(\q)$; this insight is used by \cite{francca2020conformal} to shed light on momentum-based optimization procedures such as Polyak's heavy ball method, Nesterov's accelerated-gradient method, and other alternatives. On a different note, \cite{chen2014stochastic} incorporate a friction term in their proposal of stochastic-gradient HMC in order to establish invariance of the target distribution.

%%%%%%%%%%%%%%%%%%%%%%%
% Subtracting friction
%%%%%%%%%%%%%%%%%%%%%%%

{\bfseries Subtracting friction accrues energy.\quad} If we depart from the physical analogy of Hamiltonian dynamics and consider the dynamics of a particle with an \textit{amplifying} force,
\eq{
    \dd{t}\z_t = \Om \grad_{\!\z} H(\z_t) + \G \z_t,
    \label{eq:amp-hamflow}
}
where, in contrast to \eref{eq:conf-hamflow}, the additive term $\G \z_t$ in Eq.~\eqref{eq:amp-hamflow} leads to an accrual of energy. The key observation is that, the \textit{negative coefficient of friction} (owing to the $+\gamma \p_t$ term) amplifies the motion of a particle, and the dynamics in \eref{eq:amp-hamflow} encourages the particle to \textit{flow away} from the critical points of $U(\q)$. Unlike the conservative and conformal Hamiltonian dynamics that have intuitive analogues to hockey pucks moving on frictionless or rough surfaces (i.e., preserving or dissipating energy), respectively, the amplifying Hamiltonian system does not have a physical analogue since the physical system is accruing energy as it moves. 

For $t > 0$, let $\dn_t, \up_t: \R^{2d} \to \R^{2d}$ be the integral operators solving \eref{eq:conf-hamflow} and \eref{eq:amp-hamflow} for time $t$, respectively. The following lemma provides a formal justification for the interpretation of the dynamics in \eref{eq:amp-hamflow} and  \eref{eq:conf-hamflow} as repelling and attracting dynamics, respectively.

\begin{lemma}\label{lemma:energy}
    Let $H(\z) = H(\qp) = U(\q) + \half \p\tr\S\inv\p$, then the time evolution of the Hamiltonian for Eq.~\eqref{eq:conf-hamflow} and \eref{eq:amp-hamflow}, respectively, is given by
    \eq{
        \dd{t} H(\z_t) = 
        \begin{cases}
            +\gamma \cdot \p_t\tr\S\inv\p_t & \text{if } \z_t = \up_t(\z_0),\nn\\
            -\gamma \cdot \p_t\tr\S\inv\p_t & \text{if } \z_t = \dn_t(\z_0).
        \end{cases}
        \label{eq:energy}
    }
\end{lemma}

Notably, for $\z_t = \dn_t(\z)$, the Hamiltonian $H(\z_t)$ is a Lyapunov stable function \citep[Theorem~5.16]{sastry2013nonlinear}, and all flows $\qty{\dn_t(\z): t > 0}$ will tend to critical points $\z^* = (\q^*, \p^*)$ satisfying $\grad \log\pi(\q^*) = 0$ and $\p^*=\zerov$. In contrast, for $\z_t = \up_t(\z)$ the same critical points $\z^*$ become \textit{unstable}, i.e., the flows $\qty{\up_t(\z): t > 0}$ will tend to move away from $\z^*$. Lemma~\ref{lemma:energy} is a well-known result for Hamiltonian systems with linear dissipation, and the proof is provided in Appendix~\ref{proof:lemma:energy} for completeness.

\begin{figure}[t!]
    \begin{subfigure}{0.3333\textwidth}
        \includegraphics[width=\textwidth]{\Root/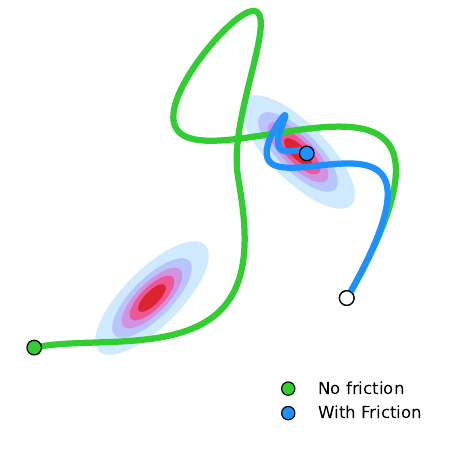}
        \caption{Attracting}
        \label{subfig:attracting}
    \end{subfigure}
    \begin{subfigure}{0.3333\textwidth}
        \includegraphics[width=\textwidth]{\Root/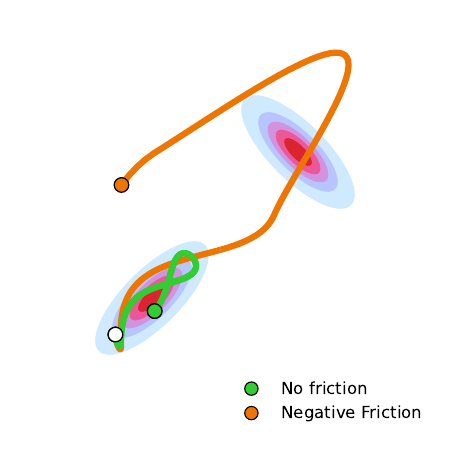}
        \caption{Repelling}
        \label{subfig:repelling}
    \end{subfigure}
    \begin{subfigure}{0.3333\textwidth}
        \includegraphics[width=\textwidth]{\Root/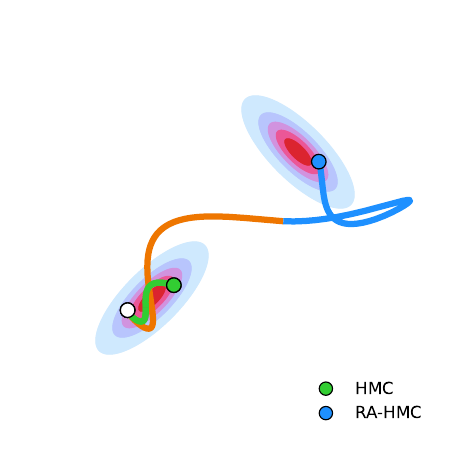}
        \caption{\haram{}}
        \label{subfig:rahmc}
    \end{subfigure}
    \caption{The dynamics of a particle   starting at $\q \in \R^2$ (denoted by $\tikzcircle[fill=white]{2.75pt}$). The first panel (a) shows a particle transition under the  conservative Hamiltonian dynamics (ending up with $\tikzcircle[fill=green]{2.75pt}$) in \eref{eq:hamflow} for fixed $\e, L$ with initial momentum $\p_1 \in \R^2$, and that under the conformal Hamiltonian dynamics (ending up with $\tikzcircle[fill=dodgerblue]{2.75pt}$) in \eref{eq:conf-hamflow} with $\gamma > 0$ and the same parameters $\e, L$. The particle under the conformal Hamiltonian dynamics with positive friction tends to be attracted to a nearby mode, unlike the one under the conservative Hamiltonian dynamics (no friction). The middle panel (b) displays the conservative Hamiltonian dynamics ($\tikzcircle[fill=green]{2.75pt}$) and amplifying Hamiltonian dynamics ($\tikzcircle[fill=orange]{2.75pt}$) of \eref{eq:amp-hamflow} for fixed $\e, L$ and $\p_2 \in \R^2$. Unlike the conservative Hamiltonian dynamics (no friction), the amplifying Hamiltonian dynamics (with negative friction) encourages a mode-repelling behavior of the particle. The last panel (c) exhibits a particle transition of \haram{} that combines the repelling and attracting dynamics ($\tikzcircle[fill=orange]{2.75pt}\!+\!\tikzcircle[fill=dodgerblue]{2.75pt}$), first being repelled from a local mode and then attracted to another mode.}
    \label{fig:friction}
\end{figure}

\begin{example}[Bivariate Anisotropic Gaussian Mixture]\label{ex:anisotropic}

Let $\pi(\q)$ be a bivariate Gaussian mixture given by
$$
\pi(\q)= 0.5 \N_2(\q\mid \muv, \Sigma_1) + 0.5 \N_2(\q\mid -\muv, \Sigma_2),
$$
where $\muv=(3, 3)^\top$, and
$$
\Sigma_1=\mat{1.0 & 0.5\\ 0.5 & 1.0}, \qq{and} \Sigma_2=\mat{1.0 & -0.5\\ -0.5 & 1.0}.
$$

Figure~\ref{subfig:attracting} contrasts the conservative Hamiltonian dynamics (in green), from the conformal Hamiltonian dynamics of~ \eref{eq:conf-hamflow} (in blue). For two trajectories starting at the same initial position, the conformal Hamiltonian dynamics in \eref{eq:conf-hamflow} encourages the particle to flow towards the isolated modes of $\pi$, whereas the conservative Hamiltonian dynamics in \eref{eq:hamflow} traverses the iso-energy contours of $\pi$. 

In contrast, flipping the sign of the friction parameter $\gamma$ transforms the \textit{attracting} states into \textit{repelling} states, and the trajectory of the particle under the amplifying Hamiltonian dynamics in \eref{eq:amp-hamflow} (in orange) traverses up the energy level-sets, as shown in Figure~\ref{subfig:repelling}. 
\end{example}

%%%%%%%%%%%%%%%%%%%%%%%%%%%%%%%%%%%%%%%C%%%%%%%%%%%%%%%%%%%%%%%%%%%%%%%%%%%%%%%%%%%%%
%%%%%% SECTION: Haram

\subsection{The proposed repelling-attracting Hamiltonian  Monte Carlo algorithm}
\label{sec:raHMCdetail}

Although conformal Hamiltonian systems present an interesting perspective for studying the modalities of the target distribution, neither of them, by themselves, are suitable for sampling from multimodal distributions since they no longer retain the key properties, \ref{p:energy}--\ref{p:symplectic}, which make  the original Hamiltonian Monte Carlo successful in practice. 

The \haram{}  algorithm combines the repelling dynamics in \eref{eq:amp-hamflow} with the attracting dynamics in \eref{eq:conf-hamflow} to propose the next state of the Markov chain, which is followed by a Metropolis-Hastings filter. Specifically, for a time  $T > 0$ and $\z_t = (\qp[t])$, the \haram{} dynamics is given by:
\begin{subequations}
    \begin{numcases}{\dd{t}\z_t = }
        \Om \grad_{\!\z} H(\z_t) + \G \z_t, & $0 < t \le \f T2$, \label{eq:haram-a}\\[5pt]
        \Om \grad_{\!\z} H(\z_t) - \G \z_t, & $\f T2 < t \le T$. \label{eq:haram-b}
    \end{numcases}
    \label{eq:haram}
\end{subequations}

For $\z = \z_0$ at time $t=0$, the position of the particle at time $t=T$ in the extended phase space is given by 
$$
\z_T = \F_T := \dn_{T/2}\circ\up_{T/2}(\z),
$$ 
and, the final \haram{} proposal, $\tz_T = \flip(\z_T)$, after a momentum flip operation, is then accepted/rejected with probability 
%\textcolor{blue}{(Tak: maybe $\up$ and $\dn$ above are swapped.)}
\eq{
    \alpha( \tz_T \mid \z ) = \min\qty\Big{1, \exp\qty\Big( H(\tz_T) - H(\z) )}.
}

In the first half of the trajectory, the particle is driven away from its current mode by implementing the mode-repelling dynamics of \eref{eq:haram-a}. Following this, in the second half of the trajectory, the particle is pulled towards a potentially new mode by the mode-attracting dynamics, as per \eref{eq:haram-b}. 
This sequence of repulsion and attraction allows for better exploration of the modes of the target distribution. For the same target distribution from Example~\ref{ex:anisotropic}, Figure~\ref{subfig:rahmc} illustrates the \haram{} dynamics, where the particle is first repelled from the current mode and then attracted to another mode.

\begin{remark}\label{rem:modes} We make some key observations about the \haram{} algorithm.
\begin{enumerate}[label=\textup{(\roman*)}, align=left, leftmargin=!, itemindent=0.5em]
    \item The Repelling-Attracting Metropolis (RAM) algorithm proposed by \cite{tak2018ram} shares a similar philosophy of combining repelling and attracting dynamics to sample from multimodal distributions. However, unlike the RAM algorithm, which is the repelling-attracting variant of random walk Metropolis-Hastings algorithm, \haram{} is based on Hamiltonian Monte Carlo. The key benefits of the proposed approach are: (a)  eliminating the need for forced-transitions via repeated sampling in low-density regions \citep[Table~1, Step~2]{tak2018ram}, and (b) better adapting to the geometry of $\supp(\pi)$, and thereby ensuring more efficient sampling from high-dimensional multimodal distributions. For detailed comparisions and results, please refer to the experiments in Section~\ref{sec:xperiments}.

    \item In light of Lemma~\ref{lemma:energy}, the \haram{} dynamics, promotes transitions among a broader class of ``modes" compared to those typically discussed in the MCMC literature. In $\R^d$ and for $k \le d$, the ``generalized modes" of $\pi$ could be interpreted as $k$-dimensional submanifolds of $\supp(\pi)$ where $\grad\log\pi(\q) = 0$. Several existing multimodal MCMC samplers (e.g., \citealp{wormhole2014,pompe2020framework}), are constructed to simulate transitions between isolated points in $\supp(\pi)$, where $\log\pi(\q)$ reaches a local maximum. These points correspond to the $0$-dimensional modes of $\pi$. In contrast, according to Lemma~\ref{lemma:energy}, \haram{} facilitates transitions between these generalized modes of $\pi$, which includes, in particular, the $0$-dimensional modes. Refer to Section~\ref{exp:low} for an instance where \haram{} efficiently samples from $\pi$ when the modes are are $1$- and $2$-dimensional submanifolds of the ambient space.
\end{enumerate}

\end{remark}

The numerical implementation of the scheme in \eref{eq:haram} is based on the \emph{conformal leapfrog integrator} (\citealp[Algorithm~3.1]{modin2011geometric} and \citealp[Eq.~3.3]{francca2020conformal}).  For a step-size $\epsilon \approx dt$, the conformal leapfrog integrator for \eref{eq:haram-a} is given by the Strang splitting
\eq{\label{eq:haram-strang}
\up_{\epsilon} = \xi^{\G}_{\epsilon/2} \circ \P^{H}_{\epsilon} \circ \xi^{\G}_{\epsilon/2},
}
where $\P^H_\epsilon$ is the usual leapfrog step for the conservative Hamiltonian dynamics from \eref{eq:xis}, and $\xi^{\G}_{\epsilon/2}(\z) = \xi^{\G}_{\epsilon/2}(\qp) = (\q,~ e^{\gamma\epsilon/2}\p)$ is the exact flow associated with the vector field $\d\z_t/\d t = \G\z_t$. Mutatis mutandis, the integrator for \eref{eq:haram-b} is obtained by negating the sign of $\G$, and is given by
$$
\dn_{\epsilon} = \xi^{-\G}_{\epsilon/2} \circ \P^{H}_{\epsilon} \circ \xi^{-\G}_{\epsilon/2}.
$$
At time $T = {L/\epsilon}$, the position of $\z_T$ in the extended phase-space is given by applying $\up_\epsilon$ repeatedly $\floor{L/2}$ times followed by $\dn_\epsilon$ repeatedly $\floor{L/2}$ times, i.e.,
\eq{
    \z_T = \F_{\epsilon, T}(\z) = \dn_{\epsilon, L/2} \circ \up_{\epsilon, L/2}(\z) = \qty(\dn_{\epsilon})^\floor{L/2} \circ \qty(\up_{\epsilon})^\floor{L/2}(\zv).
}
An algorithmic description of the \haram{} transition is provided in Algorithm~\ref{alg:haram}.

\begin{remark}
    The objective of Algorithm~\ref{alg:haram} is to first repel the particle from the current mode and then attract it to another mode, thereby encouraging better sampling from multimodal distributions. At time $T/2$, the dynamics of the particle transitions from mode-repelling to mode-attracting behavior. Mode exploration can be further encouraged by gently nudging the particle from backtracking to the same mode using the reflection principle \citep[Eq.~12]{neal2003slice}. In particular, at time $T/2$ (i.e., after Line~7 of Algorithm~\ref{alg:haram}), the transformation
    \eq{
        \mathsf{R}^{\perp}(\qp) = \qty(\q,~~ \p - 2{\grad U(\q)} \inner[\bigg]{\p, \f{\grad U(\q)}{\norm{\grad U(\q)}^2}}),\nn%\label{eq:reflection}
    }
    reflects the momentum $\p$ across the hyperplane orthogonal to the gradient of the potential energy $U(\q)$, akin to a billiard ball bouncing off a wall. In our experiments, we observed that even without the reflection principle, \haram{} is able to explore multiple modes of the target distribution efficiently, and the reflection principle is not necessary for the \haram{} algorithm to work.
\end{remark}

\begin{algorithm}[t!]
\caption{The Repelling-Attracting Hamiltonian Monte Carlo (\haram{}) algorithm.}
\label{alg:haram}
\SetAlgoLined
\SetKwFunction{Union}{Union}\SetKwFunction{haramstep}{{\normalfont \scshape conformal-Leapfrog}}
\BlankLine
\Input{Potential $U := -\log\pi$, mass matrix $\S \succ \zerov$, friction $\gamma > 0$, step size $\e > 0$, length $L >0$, sample size $N \in \mathbb{Z}_+$.}
\BlankLine
\Output{Samples $\qty{\q_1,, \q_2, \dots, \q_N} \subset \R^d$}
\BlankLine
\BlankLine
\textbf{Initialize} $\qo \in \R^d$ and \textbf{define} $H(\q, \p) := U(\q) + \f{1}{2}\p\tr\S\inv\p$\;
\BlankLine
\everypar={\nl}
% \Indp
\For{$n$ = $1 \ \KwTo \ N$}{
    Sample $\p_{n-1} \sim \N(\zerov, \Sigma)$\;
    Set $(\q,\p) \gets (\q_{n-1},\p_{n-1})$\;
    \For{$i= 1 \ \KwTo \  \floor{L/2}$}{$(\q,\p) \gets \haramstep(\q,\p,\S,\e,-\gamma)$\;}
    % Set $(\q, \p) \gets (\q, -\p)$\;
    \For{$i= 1 \ \KwTo \  \floor{L/2}$}{$(\q,\p) \gets \haramstep(\q,\p,\S,\e,+\gamma)$\;}
    Set $(\q, \p) \gets (\q, -\p)$\;
    Set $\alpha \gets \min\qty\big{1, \exp\qty\big( H(\q_{n-1}, \p_{n-1}) - H(\q,\p) )}$\;
    Sample $\uv \sim \text{Unif}(0,1)$\;
    \BlankLine
    \eIf{$\uv<\alpha$}{
        $\q_n \gets \q$\;
    }{
        $\q_n \gets \q_{n-1}$\;
    }
}
% \Indm
\textbf{return} $\qty{\q_1, \q_2, \dots, \q_N}$\;
\BlankLine
\BlankLine
\BlankLine
\SetKwProg{Fn}{Function}{:}{}
\Fn{$\haramstep(\q,\p,\S,\e,\gamma)$}{
    $\widetilde{\p} \  \gets e^{-\gamma\e/2}\p - \f{\e}{2}\grad U(\q)$\;
    $\widetilde{\q} \ \gets \q + \e\Sigma\inv \widetilde{\p}$ \;
    $\widetilde{\p} \ \gets e^{-\gamma\e/2}\qty\Big({\widetilde{\p} - \f{\e}{2}\grad U(\widetilde{\q})})$\;
    \textbf{return} $ (\widetilde{\q}, \ \widetilde{\p}) $ 
}
\end{algorithm}

% \begin{algorithm}[H]
% \caption{One step of the Conformal Leapfrog algorithm \citep[Eq.~11]{francca2020conformal}.}
% \label{alg:conformal-leapfrog}
% \SetAlgoLined
% % \LinesNumbered
% \SetKwFunction{Union}{Union}\SetKwFunction{haramstep}{{\normalfont \scshape conformal-Leapfrog}}
% \BlankLine
% \Input{Potential $U := -\log\pi$, mass matrix $\Sigma \succ \zerov$, friction $\gamma > 0$, step size $\e > 0$.}
% \BlankLine
% \Output{$(\qt, \pt)$ corresponding to one step of the conformal leapfrog integrator.}
% \BlankLine
% \BlankLine
% \SetKwProg{Fn}{Function}{:}{}
% \Fn{$\haramstep(\q,\p,\S,\e,\gamma)$}{
%     $\widetilde{\p} \  \gets e^{-\gamma\e/2}\p - \f{\e}{2}\grad U(\q)$\;
%     $\widetilde{\q} \ \gets \q + \e\Sigma\inv \widetilde{\p}$ \;
%     $\widetilde{\p} \ \gets e^{-\gamma\e/2}\qty\Big({\widetilde{\p} - \f{\e}{2}\grad U(\widetilde{\q})})$\;
%     \textbf{return} $ (\widetilde{\q}, \ \widetilde{\p}) $ 
% }
% \end{algorithm}

In order to establish that Algorithm~\ref{alg:haram} constitutes a valid Markov chain Monte Carlo algorithm, we turn our attention to investigate the properties of the \haram{} dynamics. We begin by describing a time-reveral symmetry property of \eref{eq:haram}. Following the convention in \cite{lamb1998time}, for a dynamical system $d\z_t/\d t = F(\z_t)$ with forward-time flow map $\U_t: \z_s \mapsto \z_{s+t}$, a map $R$ is defined to be a \textit{reversing symmetry} if $\U_t\qty(R(\z)) = R( \U_t\inv(\z) )$, i.e., the forward-time flow of $\U_t$ after applying $R$ is the same as $R$ applied to the backward-time flow of $\U_t$. The following result establishes a time-reversal symmetry \textit{between} $\up_t$ and $\dn_t$. 

\begin{proposition}
    Let $\dn_t, \up_t$ be the integral operators solving \eref{eq:conf-hamflow} and \eref{eq:amp-hamflow} for time $t$, respectively, and let $\flip: (\q, \p) \mapsto (\q, -\p)$ be the momentum-flip operator. Then, 
    \eq{\label{eq:time-reversal}
        \flip \circ (\up_t)^{-1} = \dn_t \circ \flip, \qq{and} 
        \flip \circ (\dn_t)^{-1} = \up_t \circ \flip \qq{for all $t \in \R$.}
    }
    In particular, for the exact \haram{} flow operator $\F_T$ and its numerical approximation $\Fel$, their compositions with the momentum-flip operator, $\flip \circ \FT \text{ and } \flip \circ \Fel$, are involutions, i.e., 
    \eq{
        \qty(\flip \circ \FT) \circ \qty(\flip \circ \FT) = \id = \qty(\flip \circ \Fel) \circ \qty(\flip \circ \Fel).\nn
    }
    \label{prop:reversibility}
\end{proposition}
\vspace*{-2em}
In other words, Proposition~\ref{prop:reversibility} states that the momentum-flip operator $\flip(\q, \p) = (\q, -\p)$ relates the forward-time flow of $\P^{\pm}_t$ with the backward-time flow of $\P^{\mp}_t$. When $\flip \circ \F_T$ is composed with itself, momentum-flip operation helps undo the forward-time flow along $\dn_T$, but leaves the momentum flipped in this intermediate state. The flow along $\dn_T$ from this intermediate state, now, reverses the forward-time flow along $\up_T$. Therefore, from \cite{tierney1994markov}, this guarantees that the resulting deterministic Markov kernel $\z \mapsto \flip \circ \F_T(\z)$ satisfies the detailed-balance condition. The proof of Proposition~\ref{prop:reversibility} is provided in Section~\ref{proof:prop:reversibility}, and holds for Hamiltonians satisfying $H(\q, -\p) = H(\q, \p)$.

Having established the reversibility of \haram{}, the following result shows that both $\F_T$ and $\F_\el$ preserve symplectic structure and volume in the extended phase space.

\begin{proposition}[Preservation of Volume and Symplectic Structure]
    Let $\omega_t$ be the symplectic $2$-form given by the wedge product ${\omega_t = \d\qt \wedge \d\pt}$. Then, \eref{eq:haram} preserves symplectic strcture only when $t=T$, i.e.,
    \eq{
        \omega_T = \omega_0,\nn
    }
    and, in particular, for all $\z \in \R^{2d}$,
    \eq{
    \jac{\FT}(\z)\tr \ \Om\inv \ \jac{\FT}(\z) = \Om\inv.\nn
    }
    Moreover, the integrators $\FT$ and $\Fel$ are volume preserving, i.e., $\abs{\det\pa{\jac{\FT}}} = \abs{\det\pa{\jac{\Fel}}} = 1$.
    \label{prop:symplecticity}
\end{proposition}

As noted in \ref{p:volume}, preservation of volume in the phase space ensures that the Metroplis-Hastings acceptance probability does not require a correction term involving the Jacobian determinant. Hence, for ${\tz_T = \flip\circ \F_T(\z)}$ or $\tz_T = \flip \circ \Fel(\z)$, the proposed state $\tz_T$ \mbox{can~be~accepted/rejected~with~probability}
$$
\alpha( \tz_T \mid \z ) = \min\qty\big{1, \exp\qty( H(\tz_T) - H(\z) )}.
$$ 

Propositions~\ref{prop:reversibility}~and~\ref{prop:symplecticity} alone are enough to guarantee that Algorithm~\ref{alg:haram} is a valid MCMC algorithm, and therefore, the samples generated by \haram{} are distributed according to the target distribution $\pi$. However, as noted in \ref{p:symplectic} and \ref{p:energy}, the access to a symplectic integrator and the energy conservation of Hamiltonian dynamics ensure high acceptance probability, even over distant trajectories. To this end, the following result establishes that the numerical scheme for \haram{}, as described in Algorithm~\ref{alg:haram}, is a second-order scheme and closely approximates the exact \haram{} dynamics.

\begin{proposition}[Approximation Error]
    For fixed $T = L\epsilon$, let $\F_T$ be the exact time-$T$ flow of the \haram{} dynamics in \eref{eq:haram}, and let $\F_\el$ be the numerical scheme in Algorithm~\ref{alg:haram}. Then, 
    \eq{
        \norm{\Fel(\z) - \F_T(\z)} = O(\e^2),\qq{and} \abs\Big{ H\qty(\Fel(\z)) - H\qty(\FT(\z)) } = O(\e^2).\nn
    }
    \vspace*{-1em}
    \label{prop:approximation}
\end{proposition}

The proof of Proposition~\ref{prop:approximation} is provided in Section~\ref{proof:prop:approximation}, and the result is similar to the case of the conservative Hamiltonian dynamics when $T = L\epsilon$ is fixed. Therefore, for fixed $T$ we don't sacrifice any numerical accuracy by departing to the conformal Hamiltonian framework. However, if $T$ depends on $\epsilon$, it can be shown via more involved and technical backward error analysis (see, e.g., \citealp{benettin1994hamiltonian,reich1999backward}), that choosing $T_\epsilon = O(\e^2 e^{C/\e})$ guarantees the same $O(\e^2)$ approximation error for $H$. This is not the case for \haram{} owing to the near-exponential growth of $H(\qt, \pt)$ in the repelling dynamics.

As far as energy conservation is concerned, the picture is more bleak; and the following result only guarantees an upper bound for the energy drift of the \haram{} dynamics.The result relies on some technical assumptions on the potential energy function $U(\q)$, and, therefore, also on the target distribution~$\pi$. 
\begin{enumerate}[label=\textup{(\textbf{A\arabic*})}]
    \item\label{a-1} The negative log density $U(\q) = -\log\pi(\q)$ is $L$-smooth, i.e., it is twice continuously differentiable and $\grad U(\q)$ is $L$-Lipschitz.
    \item\label{a-2} The Hessian $\grad^2 U(\q)$ is uniformly bounded, i.e., $\norm{\grad^2 U(\q)} \le \varpi$ for all $\q \in \R^d$.
    \item\label{a-3} $U(\q)$ is coercive, i.e., $U(\q) \to \infty$ as $\norm{\q} \rightarrow \infty$.
    \item\label{a-4} $U(\q)$ is a Morse function, and
    \item\label{a-5} $U(\q)$ satisfies the \pl{} condition, i.e., there exists a constant $\mu > 0$ such that for all critical points $\q^*$ satisfying $\grad U(\q^*) = 0$, the following inequality holds:
    $$
    \norm{\grad U(\q)}^2 \ge \mu \qty\Big(U(\q) - U(\q^*)).
    $$
\end{enumerate}

With this background, the following result provides an upper bound for the energy drift of the \haram{} dynamics.

\begin{proposition}[Bound on Energy Drift]
    Suppose $U(\q) = -\log \pi(\q)$ satisfies assumptions \ref{a-1}--\ref{a-5}. For $T > 0$, let $\F_T$  be the integral operator for the \haram{} dynamics, and let $\calo S$ be the critical set of $H(\z)$ given by
    $$
    \calo{S} := \qty{(\q^*, \zerov_d) \in \R^{2d}: \grad U(\q^*) = \zerov_d}.
    $$ 
    Then, the energy drift of the \haram{} dynamics is bounded by
    \eq{
        \abs\Big{H(\F_T(\z)) - H(\z)} \le \inf_{\uv, \vv \in \calo{S}}\abs{2H(\up_{T/2}(\z)) - H(\uv) - H(\vv)} e^{-\lambda T/2} + \sup_{\uv, \vv \in \calo S}\abs\Big{H(\uv) - H(\vv)},\nn
    }
    where $\lambda = \lambda(U, \S, \gamma) > 0$ is a constant, given in \eref{eq:lambda}, which  depends only on $U, \S$ and $\gamma$.
    \label{prop:energy}
\end{proposition}

The proof of Proposition~\ref{prop:energy} is presented in Appendix~\ref{proof:prop:energy}. While we are unable to guarantee that \haram{} conserves energy, the bound in Proposition~\ref{prop:energy} provides some control on the energy drift, and, in practice, we have observed that it is typically small. For example, Figure~\ref{fig:energy} shows the energy drift of \haram{} for random initial starting states $\z_0$ from the bivariate Gaussian mixture in Example~\ref{ex:anisotropic} for different simulation lengths $T$, and different values of the friction parameter $\gamma$. Even for long trajectories, the energy drift is small as evidenced by the $p$-values in Figure~\ref{fig:energy}. The practical implication of this is that, the optimal acceptance probability $\delta \in (0, 1)$ for \haram{} is likely to be smaller than that of HMC, which, under some strong assumptions, was shown to be $\delta \approx 0.651$ \citep{beskos2013optimal}.

\begin{remark}
    The proof of Proposition~\ref{prop:energy} is based on producing a Lyapunov function which decreases exponentially fast along the attracting stage of the \haram{} dynamics. Therefore, the assumptions \ref{a-1}--\ref{a-5} on $U(\q)$ play a key role in establishing the bound on the energy drift. We make the following remarks about these assumptions.
    \begin{enumerate}[label=\textup{(\roman*)}, align=left, leftmargin=!, itemindent=0.5em]
        \item Assumptions \ref{a-1}--\ref{a-3} are standard in the literature on dynamical systems, control theory and optimization (e.g., \citealp{sastry2013nonlinear,polyak2018optimisation}). In particular, \ref{a-1} guarantees that the trajectories are well-defined, and \ref{a-2} ensures that at each point the resulting dynamics are sufficiently well-behaved. The coercivity condition \ref{a-3} ensures that the level sets of the $H(\q, \p)$ are bounded, and, therefore, compact. For $\pi(\q)$, this requires the tails of $\pi$ to not be thicker than the exponential distribution. Similar assumptions are also needed to guarantee the ergodicity of HMC (see, e.g., \citealp{livingstone2019geometric}).
        \item Assumption~\ref{a-4} is needed to ensure that the critical points in $\calo S$ are isolated and non-singular. We note that this assumption is not restrictive, and it follows from the well-known result in Morse theory that Morse functions are dense in the space of smooth functions equipped with the uniform topology \citep[Theorem~5.27]{banyaga2004lectures}.
        \item Crucially, \ref{a-5} is needed to obtain quantitative bounds on the energy drift. The \pl{} condition ensures exponential convergence without the need for any assumptions on the convexity of $U(\q)$. This allows for a more general class of target distributions, and, in particular, for multimodal distributions. 
    \end{enumerate}
\end{remark}

\begin{figure}
    \centering
    \begin{subfigure}{0.45\textwidth}
        \includegraphics[width=\textwidth]{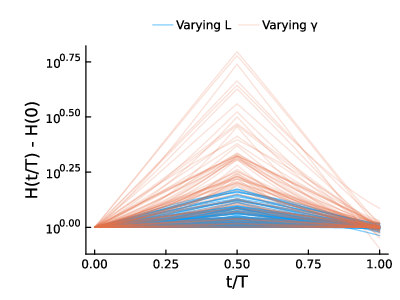}
        \caption{$H(t/T) - H(0)$ vs. $t/T$ on a $\log$ scale}
        \label{fig:energy}
    \end{subfigure}
    \quad\quad
    \begin{subfigure}{0.45\textwidth}
        \includegraphics[width=\textwidth]{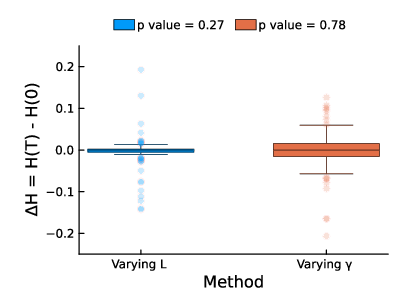}
        \caption{Boxplot of $\Delta H = H(\z_T) - H(\z_0)$}
        \label{fig:energy}
    \end{subfigure}
    \caption{The energy drift of the \haram{} dynamics for random initial starting states $\z_0$ from the bivariate Gaussian mixture in Example~\ref{ex:anisotropic}. For each of $N=100$ trajectories, the step size $\epsilon=10^{-3}$ is fixed, but the path length $L$ (blue) and friction parameter $\gamma$ (orange) are randomly sampled from $\text{Unif}\{1, 200\}$ and $\text{Unif}(0, 1)$, respectively. (Left) Plot of $H(\z_t)$ for $N=100$ trajectories where the X-axis, $t/T$, is normalized so as to depict all simulation lengths on a scale of $[0, 1]$. (Right) Boxplot of the energy drift $\Delta H = H(\z_T) - H(\z_0)$ for the same $N=100$ trajectories. The $p$-values for the  two-sided one-sample $t$-test for testing $\Delta H = 0$ vs. $\Delta H \neq 0$ shows insufficient statistical evidence that the energy drifts are statistically different from zero for randomly sampled $T$ and $\gamma$.}
    \label{fig:energy}
\end{figure}

% % \begin{remark}
% %     This result can be proven in two different ways. The first method is a direct consequence of Proposition~\ref{prop:symplecticity} above, and follows from the fact that the volume form can be represented using the symplectic form; therefore, preservation of volume follows directly from the preservation of symplectic structure. The second method is based on first principles, i.e., we explicitly characterize the distribution of $\pi(\z_t)$ for $0 \le t \le T$. Appendix~\ref{proof:prop:volume} presents both methods to prove the result.
% % \end{remark}
% %, which is presented in ,

% \begin{remark}
%     This result doesn't depend on the total time $T$ or length $L$ of the trajectory like HMC. Therefore, \haram{} retains a high acceptance probability even for long-range proposals across multiple modes.
% \end{remark}
% %, like the conservative Hamiltonian dynamics, but also a mode-jumping trajectory more frequently
% %It means that the \haram{} algorithm can generate not only a long trajectory of a proposal with high acceptance probability even in high dimensions, like the conservative Hamiltonian dynamics, but also a mode-jumping trajectory more frequently. 
% %a capability to jump between modes 

% %a proposal thatA practical meaning/interpretation of having these four properties in the context of sampling a  multimodal distribution?

% %

% %This can be discussed when we first describe the features of \haram{} in Section Section 5.3, instead of here.

\subsection{Automatic tuning procedure for \haram{} via Nesterov dual-averaging}
\label{sec:tuning}

The success of \haram{} for sampling from multimodal distributions crucially hinges on the choice of the tuning parameters $\e, \gamma, L$. Specifically,  \eref{eq:energy} guarantees that \haram{} is able to overcome high energy barriers in the multimodal target distribution $\pi$ 
%\footnote{\sout{\textcolor{red}{(Tak: Have we used this notation before? Maybe the sentence makes sense without this new notation here as you use the boldfaced mu below (not to be confused)?)}}}{\sid{$\mu$ was a typo; it's fixed now}}
by ascending iso-energy contours in the first half of the trajectory, \eref{eq:haram-a}. In particular, if $\gamma > 0$ is sufficiently large, this can be achieved with relatively small total simulation length $\epsilon L$. On the other hand, if $\gamma > 0$ is very small, then the total simulation length $\epsilon L$ must be large enough to overcome the high energy barriers. 

Even for a given choice of $\gamma$ and for an optimally tuned value of $T = \epsilon L$, the choice of $\epsilon$ is crucial for numerical stability. For instance, if the step size $\epsilon$ is too large, then, in very low density regions of the target distribution, $\log\pi(\q_{t+\epsilon}, \p_{t+\epsilon})$ can fall below machine precision and lead to numerical instability of the implementation. On the other hand, as is well-known to practicioners, if $\epsilon$ is too small then the effective computational cost of each \haram{} proposal can become prohibitively expensive.

Therefore, in order to address these concerns, we propose tuning these parameters using Nesterov dual-averaging \citep{nesterov2009dual} during the warm-up phase of the MCMC algorithm. This procedure has been shown by \cite{hoffman14a} to be effective for tuning the parameters of both HMC and the No-U-Turn Sampler (NUTS). However, unlike the setting of HMC and NUTS, the tuning of \haram{} requires the additional tuning of the friction parameter $\gamma$. 

To this end, let $\xv = (\log \epsilon,~ \log \gamma) \in \R^2$, $\gv_t(\xv) \in \R^2$ be some statistic which measures the performance of \haram{} at time $t \in \R_+$, and let $\Gv(\xv) = \bbE(\gv_t | \xv) = \lim_{T \rightarrow \infty} 1/T \sum_{t=1}^{T} \gv_t(\xv)$ denote its expected value. 
%\textcolor{blue}{distribution of $g_t$?}
%\sout{\textcolor{red}{(Tak: no dependence on t on the right-hand side?)}}\sid{It was supposed to be $F(\xv)$ and not $F_t(\xv)$. I earlier had $f$ and $F_t$ (i.e., they were interchanged), but then decided to swap to have lowercase $f$ as the stochastic quantity. I've now changed them to $g$ and $G$ because they're technically supposed to be ``subgradient'' quantities} 
Then, we can tune $\xv$ to minimize $\Gv$ using the following update rule: 
\begin{equation}
    \begin{aligned}
        \xv_{t+1} &\leftarrow \muv - \frac{\sqrt t}{\omega} \frac{1}{t_0 + t}\sum_{i=1}^n \gv_i\\
        \overline{\xv}_{t+1} &\leftarrow \eta_t \xv_{t+1} + (1 - \eta_t) \overline{\xv}_t,
    \end{aligned}
    \label{eq:nesterov}
\end{equation}
where $\overline{\xv}$ is the dual-variable associated with $\xv$, $\muv \in \R^2$ is a fixed-level that the parameters are shrunk towards, and $\omega> 0$ is a free parameter which controls the rate of shrinkage towards $\muv$, 
%the parameter
$t_0\ge 0$ controls the influence of the early iterations of the warm-up phase, and $\eta_t= t^{-k}$ for $k \in (0.5, 1]$ is a sequence of step-sizes. The main idea behind this tuning procedure is to view the statistic $\gv_t$ as the (sub-)gradient of a convex function, and then use the primal-dual averaging scheme \citep{nesterov2009dual} to tune the parameters $\xv$ so as to find the zeros of the subgradient. This, loosely, corresponds to the minimization of a convex function. As noted in \cite{hoffman14a}, if the statistic $\gv_t$ is bounded, then sequence $\overline{\xv}_t$ is guaranteed to converge to a value $\xv^*$ such that $\Gv(\xv^*) = 0$. %\sout{\textcolor{red}{(Tak: I think the notation $\bar{x}$ is not defined.)}}

Similar to the choice of the objective function for HMC in \citet[Section~3.2.3]{hoffman14a}, we propose tuning the parameters $\gamma$ and $\epsilon$ to achieve a desired Metropolis acceptance rate $\delta \in (0, 1)$. In particular, we choose
\eq{
    \gv_i = \vec{f_i, f_i} \qq{for} f_i := \delta - \min\qty\Big{1, \exp\qty\Big({ H\qty(\F_{\xv, T}(\q_i, \p)) - H(\q_i, \p) }) },\nn
}
where $T$ is a fixed simulation length, $\q_{i-1}~\in \R^d$ is the initial state, $\p \sim N(\zerov, \I_d)$ is the initial (resampled) momentum for the $i$-th iteration, and $\F_{\xv, T}(\q_i, \p)$ is the proposed state of the \haram{} trajectory in the $i$-th iteration using \eref{eq:haram} with parameters $\epsilon=\exp(\proj_1(\xv)), \gamma=\exp(\proj_2(\xv))$ and $L=T/\epsilon$. 

Under the same assumptions as \cite{hoffman14a}, i.e., if $\bbE(f_i | \xv)$ is monotone (non-increasing), it follows that $\Gv(\xv) = \bbE(\gv_i \mid \xv)$ is cyclically monotone, and can, therefore, be viewed as the (sub-)gradient of a convex function \citep[Theorem~24.8]{rockafellar1997convex}. Then, the tuning procedure in \eref{eq:nesterov} can be used to exactly achieve a desired acceptance rate $\delta$. 

In all our experiments, we set $\delta \in \qty{0.55, 0.6}$ and observed that we were able to roughly achieve the desired acceptance rate using the tuning procedure in \eref{eq:nesterov}. While the choice of $\delta = 0.65$ by \cite{hoffman14a} for HMC is motivated by the optimal acceptance rate of HMC (under fairly strong assumptions), we note that a similar choice of $\delta$ for \haram{} warrants further investigation due to the energy drift in \eref{prop:energy}.

Since there are no systematic procedures for choosing the parameters $\omega, t_0, k$ in \eref{eq:nesterov}, we use the same values as in \cite{hoffman14a} for HMC. In particular, we set $\omega = 0.05$, $t_0 = 10$, $k = 0.75$ and choose $\muv = \qty\big(\log(10\epsilon_0),~ \log(10\gamma_0))$, where $\epsilon_0$ is set using \citet[Algorithm~4]{hoffman14a} and $\gamma_0=1.0$ is set to a reasonable value which enforces traversal \textit{across} iso-energy contours of $\pi$. It is very likely that the optimal choice of these parameters for \haram{} will differ from their HMC and NUTS counterparts, but we have found them to be effective in practice, and leave the choice of these parameters as a topic for future research. 

\subsection{Extension to other variants of Hamiltonian Monte Carlo}
\label{sec:extension}

In this section, we briefly highlight how the repelling-attracting dynamics of \haram{} can be incorportated into existing extensions of HMC. \cite{lu2017relativistic} propose the Relativistic Monte Carlo algorithm, by considering the relativistic kinetic energy of the particle,
$$
K'(\p) = m_0c^2\sqrt{1 + \frac{\p\tr\p}{m_0^2c^2}},
$$
where $m_0, c > 0$ are tunable parameters which represent the rest mass and the speed of light, respectively. For the Hamiltonian $H'(\q, \p) = U(\q) + K(\p')$,
$$
\dd{t} \z_t = \Om \D_\z H'(\q_t, \p_t)
$$
represents the dynamics of a particle $\z_t$ in a relativistic setting. In practice, this amounts to replacing $\p \sim \calo{N}(\zerov, \S)$ with a symmetric hyperbolic distribution $\p \sim \pi(\p \mid m_0, c)$, and using the relativistic mass $M(\p) \in \R^{d \times d}$, the Relativistic Monte Carlo dynamics is equivalently given by
$$
\dd{t} \q_t = M(\p_t)\inv \p_t, \qq{and} \dd{t} \p_t = -\grad U(\q_t).
$$
This enables the Relativistic Monte Carlo algorithm to propose states using the leapfrog integrator. The \haram{} variant can be obtained by incorporating the repelling-attracting dynamics as follows:
$$
\dd{t} \q_t = M(\p_t)\inv \p_t, \qq{and} \dd{t} \p_t = -\grad U(\q_t) - \text{sgn}({t - T/2}) \gamma \p_t,
$$
for $0 \le t \le T$. The numerical implementation of this dynamics is obtained by including the conformal Hamiltonian steps $\xi^\Gamma_\epsilon(\q, \p) = (\q, e^{\gamma\epsilon/2}\p)$ before and after the leapfrog step in \citet[Section~2.1]{lu2017relativistic}. In a similar vein, \cite{tripuraneni2017magnetic} and \cite{brofos2020non} propose a variants of HMC whose dynamics are given by
\eq{
    \dd{t}\z_t = \underbrace{\mat{\Ev\phantom{-} & \Av \\-\Av\tr & \Gv}}_{=:\tOm(\Ev, \Gv)} \grad_{\!\z} H(\z_t),\label{eq:magnetic}
}
where $\Ev, \Gv$ are user-specified $d \times d$ skew-symmetric matrices which guarantee that $\tOm$ is a $2d \times 2d$ skew-symmetric matrix, and $\Av \in \R^{d \times d}$ is fixed so as to ensure that $\tOm$ is invertible. This differs from the canonical dynamics in \eref{eq:hamflow}, and, notably, the proposals generated by \eref{eq:magnetic} are not time-reversible in the classical sense. Instead, \eref{eq:magnetic} satisfies a weaker form of time-reversibility \citep{tripuraneni2017magnetic}, whereby negating the matrices $\Ev, \Gv$ in addition to the momentum $\p$, i.e., $\tflip(\zv, \Ev, \Gv) = (\flip(\zv), -\Ev, -\Gv)$, is needed to obtain the time-reveral dynamics. This is similar to the pesudo time-reversibility of \haram{} as noted in Proposition~\ref{prop:reversibility}, and this  additional step is necessary to guarantee that the resulting proposal is an involution, and thereby satisfy the detailed balance condition. 

\citet[Algorithm~1]{tripuraneni2017magnetic} proposes a sampler using a leapfrog-like integration scheme in the special case when $\Ev = \zerov$ and $\Gv$ is a tunable-parameter for the sampler. \citet[Algorithm~1]{brofos2020non} extends this idea to the case when, both, $\Ev$ and $\Gv$, are tunable parameters, but relies on the implicit mid-point scheme in order to generate exact proposals from $\pi$. Although the choice of $\tOm$ is crucial for the success of these samplers, the repelling-attracting dynamics of \haram{} can be incorporated into these samplers in a straightforward manner:
\begin{equation}
    \label{eq:haram-magnetic}
    \dd{t}\z_t = 
    \begin{cases}
        \tOm(\Ev, \Gv) \grad_{\!\z} H(\z_t) + \G \z_t, & 0 < t \le \f T2\\
        \tOm(\Ev, \Gv) \grad_{\!\z} H(\z_t) - \G \z_t, & \f T2 < t \le T
    \end{cases}.
\end{equation}
The exact implementation for the magnetic HMC framework can be obtained by adding an additional step before and after the leapfrog-like step in \citet[Algorithm~1] {tripuraneni2017magnetic} to account for the conformal Hamiltonian dynamics due to the $\pm \Gamma\zv_t$ term in \eref{eq:haram-magnetic}. A similar extension can be made to the non-canonical HMC framework of \citet[Algorithm~1]{brofos2020non}, and the exact details are beyond the scope of the current work. However, the key observation is that the repelling-attracting dynamics of \haram{} can be incorporated into these samplers in a straightforward manner, and can be made into valid MCMC algorithms. 

%% file: inputs/expt.tex
% !TEX root = %WORKSPACE_FOLDER%/plain.tex

\section{Numerical Illustrations}
\label{sec:xperiments}

In this section we examine the effectiveness of \haram{} in sampling from a variety of target distributions with multiple modes and complex geometry. We begin by describing the experimental setup and the target distributions used in the experiments.

\textbf{Overview.}\quad In Section~\ref{exp:da} we examine the effectiveness of the dual-averaging tuning scheme for \haram{} for a multivariate Gaussian mixture in varying dimensions. In Section~\ref{exp:low}, we compare the performance of \haram{} with four other MCMC methods on low-dimensional multimodal distributions with complex multimodal geometry, and in Section~\ref{exp:anisotropic}, we compare the same methods on high-dimensional target distributions. Lastly, in Section~\ref{exp:unimodal}, we examine the advantage of \haram{} over standard HMC for sampling from unimodal distributions. 

\textbf{Methods.}\quad We compare the performance of \haram{} to four other MCMC methods. We use HMC with dual-averaging \citep[Algorithm~5]{hoffman14a} as the baseline to demonstrate improvement, if any. We use the Repelling-Attracting Metropolis algorithm \citep[RAM,][]{tak2018ram} as a multimodal improvement of the random-walk Metropolis sampler, which requires no gradient information. In addition, we compare \haram{} to two improvements of HMC for multimodal distributions. The Wormhole HMC algorithm \citep[WHMC,][]{wormhole2014} is used as an example of an HMC extension which stores the history of the chain to enable better mode exploration, and pseudo-extended HMC \citep[PEHMC,][]{nemeth2019pseudo} is used as an example of a tempered HMC method which enables better mode transition in an extended phase-space $\R^{2dk}$ of the original target distribution. In all cases, we run each MCMC method for $5,000$ iterations after $1,000$ iterations of warm-up.

\textbf{Metrics.}\quad Throughout the numerical experiments, we evaluate the performance of each MCMC method based on three metrics: (i) the Sinkhorn distance, (ii) the average acceptance rate, and (iii) the CPU time in seconds. The average acceptance rate and the CPU time are used to measure the sampling and computational efficiency of the MCMC methods. The Sinkhorn distance is an entropic-regularized version of the Wasserstein distance \citep{cuturi2013sinkhorn} and is used as a quality measure for the resulting samples. In particular, we generate samples $\qty{\q_1, \q_2, \dots, \q_n}$ from each of the MCMC methods, and, since we know how to sample from the exact target distribution for the examples here, we generate samples $\qty{\q_1', \q_2', \dots, \q_n'} \sim \pi$. We evaluate the Sinkhorn distance between their empirical measures, $\w(\nu_n, \nu'_n)$, as a measure of sample quality. While this introduces some additional error from the empirical measure, $\nu'_n$, we have found it to be a good metric for examining sample quality from multimodal distributions since the Wasserstein distance is intimately connected to the geometry of the target distribution, unlike other alternatives, e.g., the kernel Stein discrepancy. See, e.g., \cite{panaretos2019statistical} for an overview of the Wasserstein distance and its applications in statistics. In only the last example, where the target is a unimodal Gaussian distribution, we additionally compare auto-correlation functions (ACF). We do not compare ACFs in the multimodal examples since they can be misleading in distingushing between Markov chains that are stuck in a local modes \citep{agarwal2022globally}.

\textbf{Hyperparameter tuning.}\quad In order to ensure that each method is properly tuned, we employ a number of techniques. For \haram{}, we use the procedure outlined in Section~\ref{sec:tuning} to tune the parameters $\e$, $\gamma$, and $L$ with $\delta \in \qty{0.55, 0.6}$. For HMC, we use a similar dual-averaging approach to tune $\e$ and $L$, with $\delta$ fixed at $0.6$ \citep{hoffman14a}. For both \haram{} and HMC we select the expected simulation length, $\lambda = L\epsilon$ to ensure an acceptance probability between $40\%$ and $80\%$. For the RAM algorithm, we use a Gaussian distribution with isotropic covariance $\sigma^2\I_d$ as the proposal kernel, and we tune the value of $\sigma^2$ so that the acceptance rate is less than $90\%$. In some cases---particularly in higher dimensions---we report results of RAM with very low acceptance rates so as to avoid sticky samples. For PEHMC, similar to the experiments in \cite{nemeth2019pseudo}, we set the number of pseudo-samples, $N$, to be between $50$ and $200$, and we tune the leapfrog parameters $\epsilon$ and $L$ in order to obtain an acceptance rate greater than $40\%$. Finally, to tune WHMC we set the number of fixed-point steps to $10$ and fix the \textit{world distance} parameter and the {influence factor} for the mollifier at 1.0. We use an L-BFGS search with $100,000$ max iterations for all {regeneration routines}. We also tune the leapfrog parameters $\epsilon$ and $L$, as well as the temperature $T$ for the {tempered residual potential energy} of WHMC, to achieve an average acceptance rate greater than 40\%.

%%%%%%%%%%%%%%%%%%%%%%%%%%%%%%%%%%%%%%%%%%%%%%%%
%%%%%% Dual Averaging
%%%%%%%%%%%%%%%%%%%%%%%%%%%%%%%%%%%%%%%%%%%%%%%%

\begin{figure}[t!]
    \includegraphics[width=\textwidth]{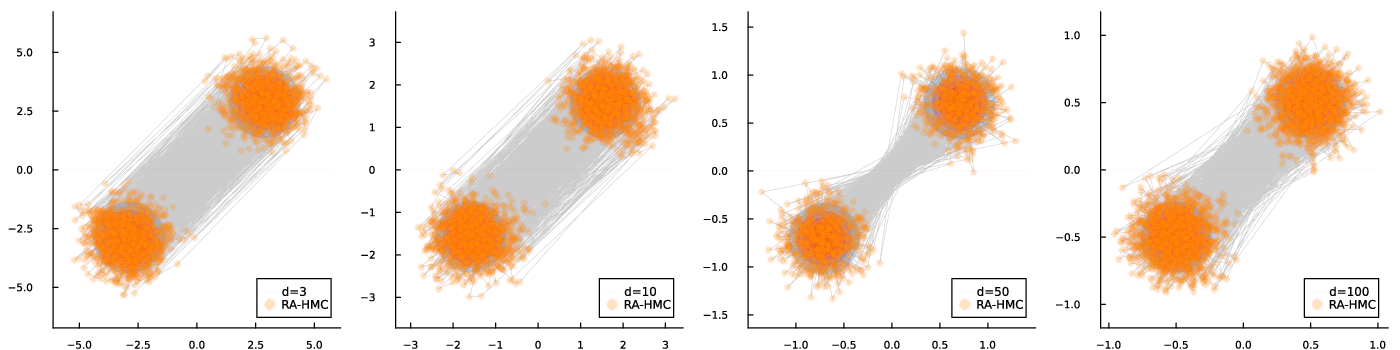}
    \caption{Results of parameter tuning for the \haram{} algorithm using dual averaging over $n=10,000$ iterations for $d \in \{3, 10, 50, 100\}$. The scatterplots show the generated samples with gray lines indicating traces between consecutive samples.}
    \label{fig:tuning}
\end{figure}

\begin{table}[t!]\centering
    \caption{\haram{} with dual averaging ($\delta=0.65$) for bimodal, isotropic Gaussian mixture. }\label{table:dual-averaging}
    \small
    \begin{tabular}{l|rrrrr}\toprule
        Dimension &$d=3$ &$d=10$ &$d=50$ &$d=100$ \\\midrule
        $\mathbf{W_2}$ \textbf{metric}  &0.16 &0.12 &0.27 &0.17 \\
        \textbf{Acceptance Rate} &63.9\% &65.1\% &36.7\% &74.9\%\\
        \textbf{CPU Time} & $27.9$s &$19.4$s &$1861.3$s &$7996.64$s\\
        \bottomrule
    \end{tabular}
\end{table}

\subsection{Performance of dual averaging for {\scriptsize RA}HMC in high dimensions}
\label{exp:da}

In this experiment, we evaluate the performance of auto-tuning the \haram{} parameters $\e, \gamma, L$ via dual averaging. The target distribution is a bimodal mixture of Gaussian distributions,
$$
\pi(\q) \propto \N_d(\xv \mid \muv_1,~ \sigma^2 \I_d) + \N_d(\xv \mid \muv_2,~ \sigma^2 \I_d),
$$
where $\muv_1 = -\muv_2 = (5 / \sqrt{d})\onev_d$ and the variance $\sigma^2 = 1/d$ are rescaled by the dimension. From Proposition~\ref{prop:mode-transition}, this effectively isolates the effect of increasing the dimension. For this example, the $2$-Wasserstein metric between the target distribution $\pi$ and any of its individual mixture components is a fixed value,
\eq{
    \w\qty(\N_d(\muv_1,~ \sigma^2 \I_d),~ \pi) = \w\qty(\N_d(\muv_2,~ \sigma^2 \I_d),~ \pi) = \frac{1}{\sqrt 2} \norm{\muv_1 - \muv_2}_2 \approx 7.07.\nn
}
In particular, if a chain is stuck at a mode, the $\w$ metric is expected to be close to $7.07$. The results from the auto-tuning \haram{} via dual-averaging are summarized in Table~\ref{table:dual-averaging} and Figure~\ref{fig:tuning}. The results for the $\w$ metric in Table~\ref{table:dual-averaging} are well below the threshold of $7.07$, indicating that \haram{} discovers the two modes and that the samples are well-mixed even in $\R^{100}$. In Figure~\ref{fig:tuning}, the two-dimensional projections of the samples show that the sample \haram{} successfully makes frequent mode-jumping transitions.

%%%%%%%%%%%%%%%%%%%%%%%%%%%%%%%%%%%%%%%%%%%%%%%%
%%%%%% Low dimensions
%%%%%%%%%%%%%%%%%%%%%%%%%%%%%%%%%%%%%%%%%%%%%%%%

\subsection{Low dimensional examples}
\label{exp:low}

\begin{figure}[t!]
    \begin{subfigure}{0.245\textwidth}
        \includegraphics[width=\textwidth]{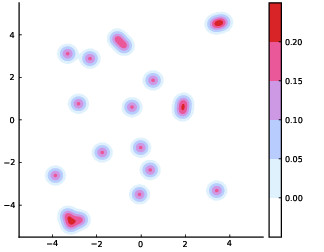}
        \caption{Benchmark dataset}
        \label{subfig:benchmark}
    \end{subfigure}
    \begin{subfigure}{0.245\textwidth}
        \includegraphics[width=\textwidth]{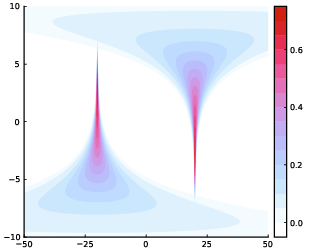}
        \caption{Multimodal Funnel}
        \label{subfig:funnel}
    \end{subfigure}
    \begin{subfigure}{0.245\textwidth}
        \includegraphics[width=\textwidth]{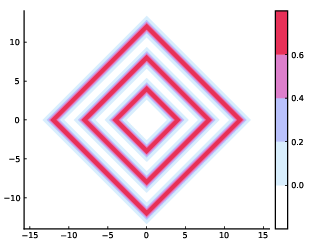}
        \caption{Concentric $\ell_1$ balls ($\R^2$)}
        \label{subfig:concentric}
    \end{subfigure}
    \begin{subfigure}{0.245\textwidth}
        \includegraphics[width=\textwidth]{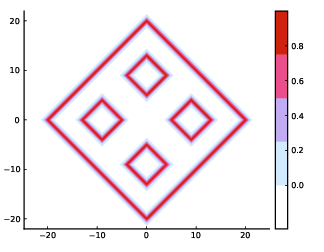}
        \caption{Nested $\ell_1$ balls ($\R^2$)}
        \label{subfig:nested}
    \end{subfigure}
    \caption{Ground truth for the low-dimensional target distributions in Section~\ref{exp:low}}
    \label{fig:low}
\end{figure}

\begin{table}[t!]
    \centering
    \caption{Results for Low dimensional examples from Section~\ref{exp:low}. Since the RAM algorithm does not make use of gradients, the information about the number of gradient evaluations $L$ and that about the CPU time per each gradient evaluation are not available (denoted by NA).  }\label{table:low}
    \small
    \begin{tabular}{cl|r C{6em} C{8em}  R{7em} R{6em} r}\toprule
    \textbf{Target} 
    &\textbf{Method} 
    &\textbf{$\mathbf{W_2}$ metric} 
    &\textbf{Acceptance Rate} & \textbf{\# Gradients per step ($L$)}
    &\specialcell{\bfseries CPU Time \\ \bfseries per iteration \\
    \bfseries (Time /$L$)} \\
    %&\specialcell{\bfseries Memory \\ (per $L$)} \\ 
    \midrule
        \multirow{5}{*}{\textbf{Benchmark}} &\textbf{HMC} &0.337 &92.1\% &27 &24.2s (0.9s) %&\textbf{2.97gb (0.11gb)}
        \\
    &\textbf{\haramalt{}} &\textbf{0.011} &57.4\% &239 &95.6s (0.4s) 
    %&0.14gb 
    \\
    &\textbf{RAM*} &0.025 &22.1\% &NA & \textbf{1.8s  (NA)}
    %&0.35gb 
    \\
    &\textbf{PEHMC} &0.152 &56.7\% &7 &938.0s (134s)
    %&20.86gb 
    \\
    &\textbf{WHMC} &0.051 &84.3\% &7 &139.3s (19.9s)
    %&10.49gb 
    \\ 
    
    \midrule
    \multirow{5}{*}{\specialcell{\textbf{Multimodal} \\ \textbf{Funnel} }} &\textbf{HMC} &26.823 &42.7\% &55 &8.15s (0.15s)
    %&\textbf{0.05gb} 
    \\
    &\textbf{\haramalt{}} &\textbf{1.519} &22.2\% &5471 &1094.00s (0.20s) 
    %&0.05gb 
    \\
    &\textbf{RAM*} &4.290 &5.3\% &NA &\textbf{0.76s (NA)}
    %&0.27gb 
    \\
    % &\textbf{PEHMC} &0.152 &56.7\% &7 &134s &20.86gb \\
    &\textbf{WHMC} &27.383 &94.9\% &100 &117.00s (1.17s)
    %&0.28gb 
    \\ 
    \midrule
    \multirow{5}{*}{\specialcell{\bfseries Concentric \\ ($d=2$)}} &\textbf{HMC} &4.295 &29.4\% &27 &10.8s (0.4s)
    %&0.19gb 
    \\
    &\textbf{\haramalt{}} &\textbf{0.191} &63.2\% &100 & 90.0s (0.9s)
    %&0.24gb 
    \\
    &\textbf{RAM*} &0.407 &43.8\% &NA &\textbf{1.1s (NA)} %
    %&\textbf{0.14gb} 
    \\
    &\textbf{PEHMC} &4.418 &84.6\% &20 &64.0s (3.2s)
    %&0.87gb 
    \\
    &\textbf{WHMC} &4.233 &60.0\% &10 &24.0s (2.4s)
    %&0.50gb 
    \\ 
    \midrule
    \multirow{5}{*}{\specialcell{\bfseries Concentric \\ ($d=3$)}} &\textbf{HMC} &3.781 &67.3\% &45 & 22.5s (0.5s)
    %&0.12gb 
    \\
    &\textbf{\haramalt{}} &\textbf{0.647} &52.8\% &510 & 459.0s (0.9s)
    %&0.24gb 
    \\
    &\textbf{RAM*} &1.618 &27.9\% &NA &\textbf{0.1s (NA)}
    %&\textbf{0.05gb} 
    \\
    &\textbf{PEHMC} &21.640 &79.4\% &7 &12.6s (1.8s)
    %&0.63gb 
    \\
    &\textbf{WHMC} &19.696 &86.1\% &20 &30.0s (1.5s)
    %&0.39gb 
    \\ 
    \midrule
    \multirow{5}{*}{\specialcell{\bfseries Nested \\ ($d=2$)}} &\textbf{HMC} &14.088 &90.3\% &12 & 1.2s (0.1s) 
    %&\textbf{0.05gb} 
    \\
    &\textbf{\haramalt{}} &\textbf{0.272} &47.3\% &61 &12.2s (0.2s)
    %&0.05gb 
    \\
    &\textbf{RAM*} &0.372 &46.5\% &NA &\textbf{1.4s (NA)}
    %&0.15gb 
    \\
    &\textbf{PEHMC} &1.036 &50.0\% &20 &150s (7.5s)
    %&1.59gb 
    \\
    &\textbf{WHMC} &13.742 &92.0\% &10 &24.0s (2.4s)%
    %&0.32gb 
    \\ 
    \midrule
    \multirow{5}{*}{\specialcell{\bfseries Nested \\ ($d=3$)}} &\textbf{HMC} &3.868 &65.0\% &289 & 28.9s (0.1s)
    %&\textbf{0.02gb} 
    \\
    &\textbf{\haramalt{}} &\textbf{1.159} &63.7\% &1655 & 165.5s (0.1s)
    %&0.04gb 
    \\
    &\textbf{RAM*} &2.089 &32.9\% &NA &\textbf{0.1s (NA)}
    %&0.07gb 
    \\
    &\textbf{PEHMC} &64.030 &80.8\% &20 &120.0s (6.0s) %&51.52gb 
    \\
    &\textbf{WHMC} &97.863 &79.0\% &25 &25.0s (1.0s) %&0.21gb 
    \\
    \bottomrule
    \end{tabular}
    \newline
    \scriptsize{*Since the RAM algorithm is gradient-free, the total CPU time and total memory are reported.}
\end{table}

In this experiment, we compare the performance of each method on low-dimensional distributions with complex multimodal geometry, as shown in Figure~\ref{fig:low}. We consider the following~four~target~distributions.\linebreak

\begin{description}[leftmargin=*, itemindent=0em]
    \item[\textbf{(A) Benchmark Dataset.}]\label{eg-low-1} Our first example is a mixture of twenty bivariate Normal distributions with equal weight as the target distribution \citep{kou2006equi}. The mixture is defined as
    $$
    \pi(\q)\propto\sum_{j=1}^{20}\exp\!\left(-\frac{(\q-\mu_j)^\top(\q-\mu_j)}{2\sigma^2}\right)
    $$
    for $\q^\top=(q_1, q_2)\in \mathbb{R}^2$ with $\sigma^2=0.01$. The mean vectors of the twenty components are available in \cite{kou2006equi}. The target distribution is plotted in Figure~\ref{subfig:benchmark}. 
    \item[\textbf{(B) Multimodal Funnel.}]\label{eg-low-2} Our second example is a multimodal generalization of Neal's funnel \citep{neal2003slice}, and the $2$-dimensional target distribution $\pi$ is given by 
    \eq{
        \pi(q_1, q_2) = \pi_1(q_1 | q_2, c, \mu) \times \pi_2(q_2 | \mu, \sigma),\nn
    }
    where $\pi(q_2 | \mu, \sigma) \propto \N(\mu - 5, \sigma) + \N(\mu + 5, \sigma)$ and $\pi(q_1 | q_2, c, \sigma) \propto \N(c, e^{c - \mu - q_2/2 }) + \N(-c, e^{c - \mu - q_2/2})$.
    % \eq{
    %     \pi(q_2 | \mu, \sigma) &\propto \N(\mu - 5, \sigma) + \N(\mu + 5, \sigma), \qq{and}%\nn\\ 
    %     \pi(q_1 | q_2, c, \sigma) &\propto \N(c, e^{c - \mu - q_2/2 }) + \N(-c, e^{c - \mu - q_2/2}).\nn
    % }
    We set $\mu=3$, $\sigma=c=1$, and Figure~\ref{subfig:funnel} displays the corresponding target distribution of ($q_1, q_2$). 
    \item[\textbf{(C) Concentric $\ell_1$ Balls.}]\label{eg-low-3} In our remaining examples, we consider target distributions with subtle geometric structures in addition to multimodality. The first of these is a mixture of distributions supported on the boundary of three \textit{concentric} $\ell_1$ balls in $\R^d$, and the target $\pi(\q)$ is given by
    \eq{
        \pi(\q) \propto \sum_{i=1}^3 \exp\qty( - (\norm{\q}_1 - r_i)^2 / 2\sigma^2 ),\nn
    }
    where $r_1 = 4, r_2=8, r_3=16$ and $\sigma=0.5$. We consider the cases when $d=2$ and $d=3$. See Figure~\ref{subfig:concentric} for the target distribution in $\R^2$.
    \item[\textbf{(D) Nested $\ell_1$ Balls.}]\label{eg-low-4} Our final example is a mixture of distributions supported on the boundary of \textit{nested} $\ell_1$ balls in $\R^d$. The target distribution, $\pi(\q)$, is given by
    \eq{
    \pi(\q) \propto \sum_{i=1}^5 \exp\qty( - (\norm{\q - \muv_i}_1 - r_i)^2 / 2\sigma^2 ),\nn
    }
    where $\muv_1 = \zerov$, $r_1 = 20$, and $\norm{\muv_i}_1 = 2, r_i = 2$ for $2 \le i \le 5$. Similar to example (C), we consider the cases when $d=2$ and $d=3$, and the target distribution when $d=2$ is shown in Figure~\ref{subfig:nested}.
\end{description}

\bigskip

\begin{figure}[t!]
    \centering
    \begin{subfigure}{0.25\textwidth}
        \caption*{Ideal}
        \includegraphics[width=\textwidth]{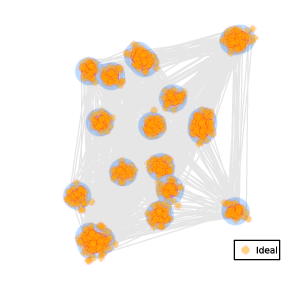}
    \end{subfigure}
    \begin{subfigure}{0.25\textwidth}
        \caption*{HMC}
        \includegraphics[width=\textwidth]{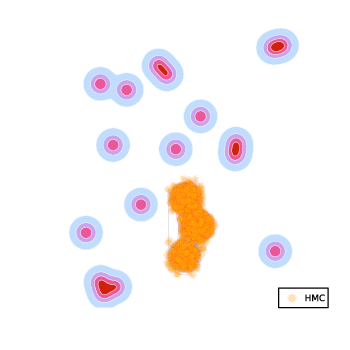}
    \end{subfigure}
    \begin{subfigure}{0.25\textwidth}
        \caption*{\haram{}}
        \includegraphics[width=\textwidth]{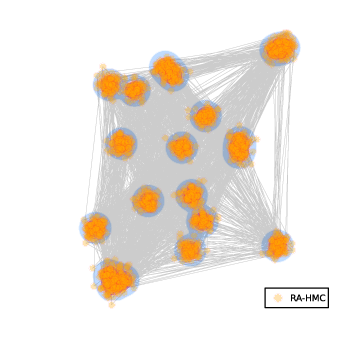}
    \end{subfigure}

    \begin{subfigure}{0.25\textwidth}
        \caption*{RAM}
        \includegraphics[width=\textwidth]{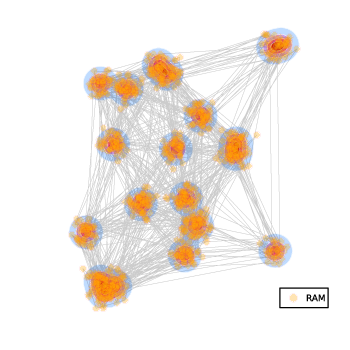}
    \end{subfigure}
    \begin{subfigure}{0.25\textwidth}
        \caption*{PEHMC}
        \includegraphics[width=\textwidth]{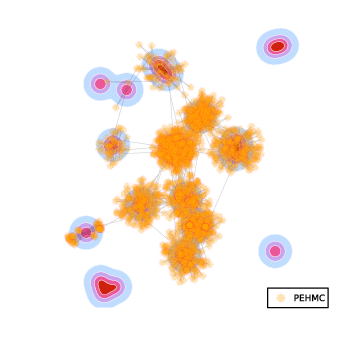}
    \end{subfigure}
    \begin{subfigure}{0.25\textwidth}
        \caption*{WHMC}
        \includegraphics[width=\textwidth]{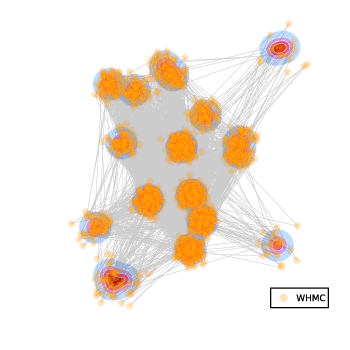}
    \end{subfigure}
    % \begin{subfigure}{0.25\textwidth}
    %     \caption*{WHMC (known modes)}
    %     \includegraphics[width=\textwidth]{\Root/sims/benchmark/whmc_known.\imageformat}
    % \end{subfigure}
    \caption{Trace plots of the last 1,000 iterations obtained by the five different MCMC samplers for the benchmark target distribution in Section~\ref{exp:low}\,(A). The plot in the top-left corner shows the ideal trace plot if the sampler was able to explore all the modes.}
    \label{fig:benchmark}
\end{figure}

The results for the experiment are summarized in Table~\ref{table:low}. The results show that \haram{} outperforms the other methods in terms of the $\w$ metric for all the examples. This, however, comes at the price of a higher computational cost, as the auto-tuning procedure for \haram{} tends to favor longer trajectories, and therefore requires more gradient evaluations per step. RAM is able to achieve reasonably low $\w$ accuracy with better CPU time, but comes at the cost of low acceptance rates. This can be attributed to the forced mode transitions with a random-walk kernel \citep[Eq.~3]{tak2018ram} which can lead to rejections when the target distribution is supported on lower dimensional structures. WHMC performs competitively when the $\pi(\q)$ has well-separated and isolated modes, as in the benchmark dataset, but performs poorly when the modes of the target distribution are, in fact, entire submanifolds as is the case for examples~(C)~and~(D).

\begin{figure}[t!]
    \centering
    \includegraphics[width=\textwidth]{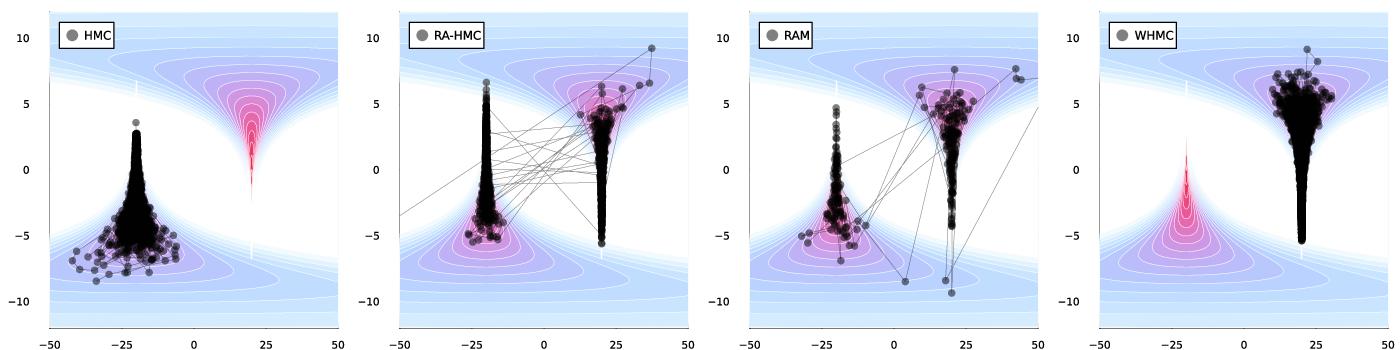}
    \caption{Scatterplots of the last 1,000 iterations obtained by four  different samplers for the multimodal funnel distribution in Section~\ref{exp:low}. \haram{} and RAM are the only samplers that move between two modes. In particular, \haram{}  samples the sharp tips of the funnels more accurately than RAM does.}
    \label{fig:multimodal-funnel}
\end{figure}

The traceplots for the benchmark dataset in example~(A) are displayed in Figure~\ref{fig:benchmark}, and show that \haram{} is able to explore the twenty modes effectively and comes closest to the ideal sampling scenario. For the multimodal funnel in example~(B), \haram{} and RAM are the only two algorithms that are able to move between the two high-density regions, but \haram{} is able to sample the sharp tips of the funnels more accurately than RAM does, as shown in Figure~\ref{fig:multimodal-funnel}.

\begin{figure}[t!]
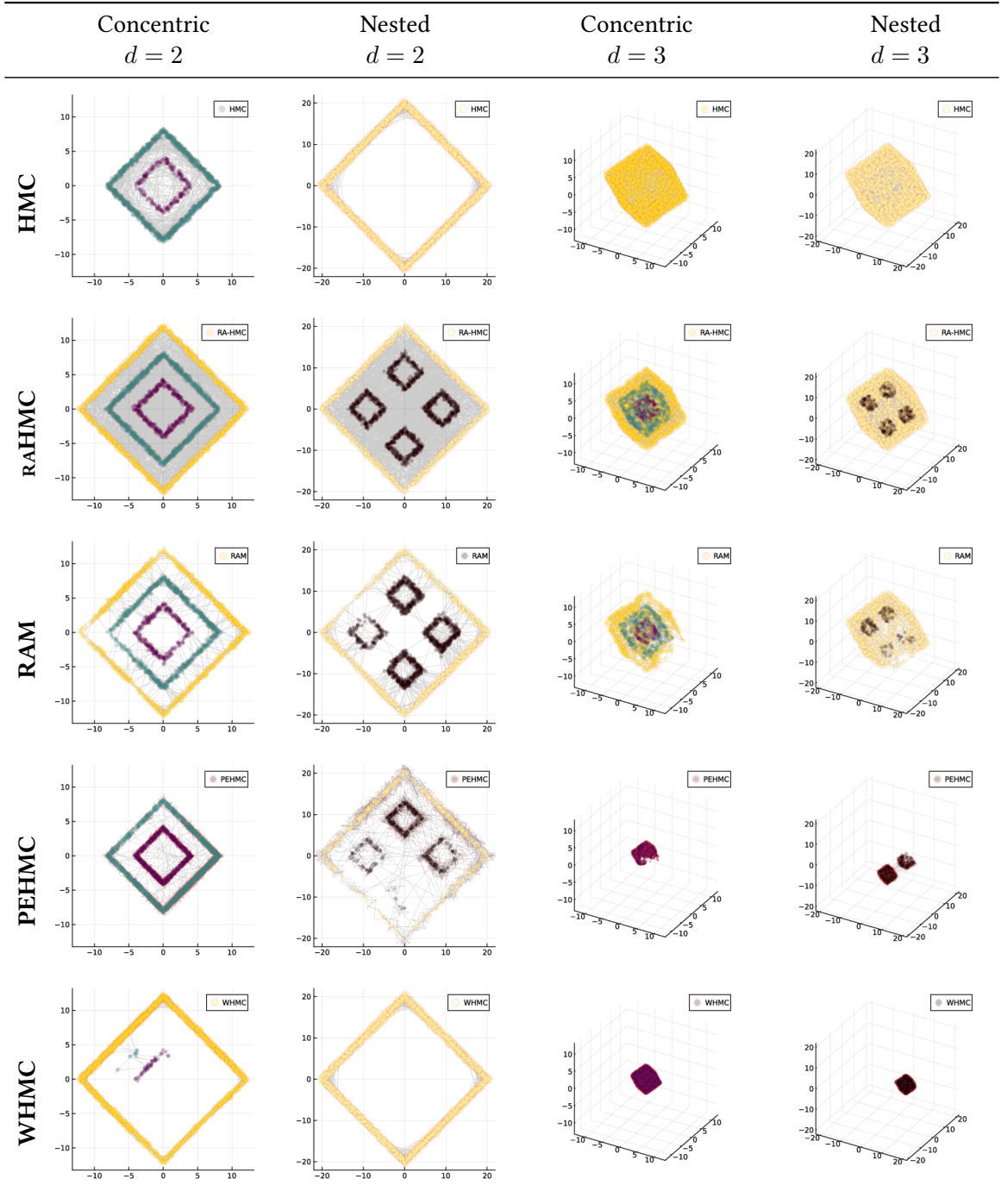

    % \begin{table}[htbp]
        \centering
        \begin{tabular}{>{\bfseries}m{0.001\textwidth} >{\centering}m{0.2\textwidth} >{\centering}m{0.2\textwidth} >{\centering}m{0.2\textwidth} >{}m{0.2\textwidth} }
            \toprule
            \ & Concentric & Nested & Concentric & \hspace*{4em}Nested \\
            \ & $d=2$ & $d=2$ & $d=3$ & \hspace*{4em}$d=3$ \\
            \midrule
            \rotatebox[origin=c]{90}{HMC} & \tabfig{0.2}{\Root/sims/circles/concentric/scatter-2d-hmc.\imageformat} & \tabfig{0.2}{\Root/sims/circles/nested/scatter-2d-hmc.\imageformat} & \tabfig{0.2}{\Root/sims/circles/concentric/scatter-3d-hmc.\imageformat} & \tabfig{0.2}{\Root/sims/circles/nested/scatter-3d-hmc.\imageformat} \\
            \rotatebox[origin=c]{90}{{\haramalt{}}} & \tabfig{0.2}{\Root/sims/circles/concentric/scatter-2d-rahmc.\imageformat} & \tabfig{0.2}{\Root/sims/circles/nested/scatter-2d-rahmc.\imageformat} & \tabfig{0.2}{\Root/sims/circles/concentric/scatter-3d-rahmc.\imageformat} & \tabfig{0.2}{\Root/sims/circles/nested/scatter-3d-rahmc.\imageformat} \\
            \rotatebox[origin=c]{90}{RAM} & \tabfig{0.2}{\Root/sims/circles/concentric/scatter-2d-ram.\imageformat} & \tabfig{0.2}{\Root/sims/circles/nested/scatter-2d-ram.\imageformat} & \tabfig{0.2}{\Root/sims/circles/concentric/scatter-3d-ram.\imageformat} & \tabfig{0.2}{\Root/sims/circles/nested/scatter-3d-ram.\imageformat} \\
            \rotatebox[origin=c]{90}{PEHMC} & \tabfig{0.2}{\Root/sims/circles/concentric/scatter-2d-pehmc.\imageformat} & \tabfig{0.2}{\Root/sims/circles/nested/scatter-2d-pehmc.\imageformat} & \tabfig{0.2}{\Root/sims/circles/concentric/scatter-3d-pehmc.\imageformat} & \tabfig{0.2}{\Root/sims/circles/nested/scatter-3d-pehmc.\imageformat} \\
            \rotatebox[origin=c]{90}{WHMC} & \tabfig{0.2}{\Root/sims/circles/concentric/scatter-2d-whmc.\imageformat} & \tabfig{0.2}{\Root/sims/circles/nested/scatter-2d-whmc.\imageformat} & \tabfig{0.2}{\Root/sims/circles/concentric/scatter-3d-whmc.\imageformat} & \tabfig{0.2}{\Root/sims/circles/nested/scatter-3d-whmc.\imageformat} \\
            \bottomrule
        \end{tabular}
        \caption{Scatterplots of the last 1,000 iterations obtained by five  different samplers for the concentric and nested $\ell_1$ balls in both $d=2$ and $d=3$ in Section~\ref{exp:low}. \haram{} and RAM are the only samplers that capture the lower dimensional submanifold structure in addition to multimodality. Although the sampling quality of RAM  deteriorates noticeably in $d=3$, that of \haram{} is stable regardless of the dimension.}
        \label{fig:l1}
    % \end{table}
    \end{figure}

The results for the concentric and nested $\ell_1$ balls in examples~(C)~and~(D) are shown in Figure~\ref{fig:l1}. While both RAM and \haram{} are competitive when $d=2$ and are able to discover and sample from the $1$-dimensional modes effectively, increasing the dimension to $d=3$ makes RAM less effective; in fact, the acceptance probability drops from $~50\%$ to $~30\%$. When $d=2$, HMC is able to move between two of the three generalized modes for the concentric $\ell_1$ case, but is unable to move between the modes for the nested $\ell_1$ case. WHMC is able to discover one of the high-density regions, but when the inner-optimization routine is run from an existing high-density region, it converges to the density ridge associated with the same $\ell_1$ ball. Therefore, it is unable to discover or move between these generalized modes. PEHMC, on the other hand, performs well when $d=2$, but when $d=3$ it is unable to discover the modes of the target distribution effectively. This may be, perhaps, attributed to the fact that PEHMC relies on the assumption that the cartesian product of the phase-space leads to intersections between different modes---which may not be the case for the nested $\ell_1$ balls. It also likely, however, that a more careful tuning of the WHMC and PEHMC hyperparameters could lead to better performance; we were unable to find such a tuning in this work.

%%%%%%%%%%%%%%%%%%%%%%%%%%%%%%%%%%%%%%%%%%%%%%%%
%%%%%% High dimensions
%%%%%%%%%%%%%%%%%%%%%%%%%%%%%%%%%%%%%%%%%%%%%%%%

\subsection{Anisotropic Gaussian Mixture in High Dimensions}
\label{exp:anisotropic}

In this experiment, we consider a generalization of the bivariate Gaussian mixture model from Example~\ref{ex:anisotropic}. For $\q \in \R^d$, the target distribution is
\eq{%\label{eq:ex1target}
    \pi(\q)= 0.5 \N_d(\q\mid \bd,~ \Sigma_1) + 0.5 \N_d(\q\mid -\bd,~ \Sigma_2),\nn
}
where $b=2$, and the covariance matrices are given by $[\Sigma_1]_{i,j} = 0.75^{\abs{i-j}}$ for $i,j = 1, 2, \dots, d$ and $\Sigma_2 = 2\I_d-\Sigma_1$. We compare the performance of four MCMC methods when $d \in \{2, 10, 20, 50, 100\}$, and the results are summarized in Table~\ref{table:anisotropic}. We omit the WHMC algorithm from this experiment as we were unable to find an appropriate tuning of the hyperparameters for high-dimensions. The covariance structure of the target distribution is particularly challenging for gradient based methods, as the principal eigenvectors of the covariance matrices are perpendicular to each other. \haram{} maintains the smallest $\w$ metric for all the dimensions except $d=2$, and the acceptance rate is also stable across the dimensions. As seen in Section~\ref{exp:low} RAM is able to achieve reasonably low $\w$ accuracy with better CPU time when $d=2$. HMC and PEHMC also perform competitively when $d=2$, but their performance deteriorates as the dimension increases. 

\begin{table}[t!]
    \caption{Summary of the sampling results obtained by four different sampling methods for the anisotropic Gaussian mixtures in Section~\ref{exp:anisotropic}. The figure on the right-hand side visualizes each method's $\w$ values over dimensions. It shows how quickly the performance of the  samplers can deteriorate when the dimension increases, while \haram{} maintains the quality of the resulting sample in a relatively stable manner.}\label{table:anisotropic}
    \vspace*{-1.5em}
    \begin{minipage}{0.67\linewidth}
    \centering
    \footnotesize
    \begin{tabular}{cl | rrrrrr }\toprule
        {} &\textbf{Method} &${d=2}$ &${d=10}$ &${d=20}$ &${d=50}$ &${d=100}$ \\\midrule
        \multirow{4}{*}{\bfseries $\mathbf W_2$ metric} &HMC &3.33 &5.55 &8.07 &13.11 &19.56 \\
        &\haram{} &0.39 &\textbf{0.77} &\textbf{1.35} &\textbf{1.99} &\textbf{3.50} \\
        &RAM &\textbf{0.26} &5.46 &8.12 &13.50 &20.21 \\
        &PEHMC &0.71 &3.79 &6.64 &12.41 &18.83 \\ 
        \midrule
        \multirow{4}{*}{\specialcell{\textbf{Acceptance} \\ \textbf{Rate}}} &HMC &71.7\% &82.8\% &59.6\% &91.1\% &97.5\% \\
        &\haram{} &49.2\% &65.7\% &62.6\% &67.2\% &71.4\% \\
        &RAM &31.9\% &11.8\% &43.7\% &29.1\% &50.1\% \\
        &PEHMC &92.0\% &63.6\% &63.5\% &88.5\% &69.2\% \\
        %\midrule
        %\multirow{4}{*}{\bfseries ESS} &HMC &\textbf{2890.54} &28.64 &\textbf{367.08} &\textbf{348.31} &\textbf{317.10} \\
        %&\haram{} &1401.32 &\textbf{777.28} &347.73 &10.96 &255.94 \\
        %&RAM &128.16 &46.17 &21.72 &14.51 &18.48 \\
        %&PEHMC &207.60 &219.65 &276.21 &126.45 &107.68 \\ 
        \bottomrule
        % \vspace*{1em}
        \end{tabular}
        % \caption{tab caption}
        % \label{table:student}
    \end{minipage}
    \begin{minipage}{0.3\textwidth}
    \centering
    \vspace*{1em}
    \includegraphics[width=1.0\textwidth]{\Root/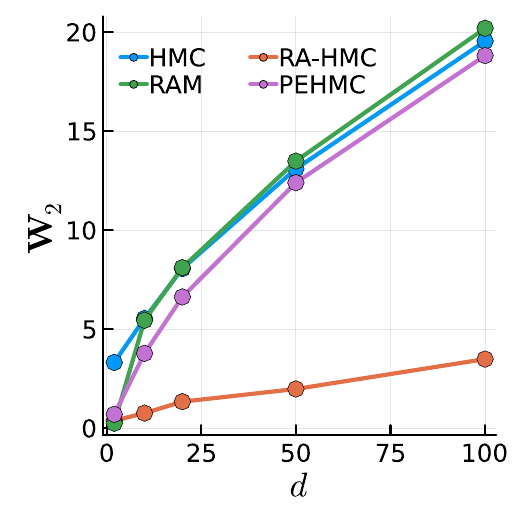}
    % \vspace*{-1em}
    % \captionof{figure}{fig caption}
    % \label{fig:image}
    \end{minipage}
    \vspace*{-1em}
\end{table}

%%%%%%%%%%%%%%%%%%%%%%%%%%%%%%%%%%%%%%%%%%%%%%%%
%%%%%% Unimodal Gaussian
%%%%%%%%%%%%%%%%%%%%%%%%%%%%%%%%%%%%%%%%%%%%%%%%

\subsection{Sampling from Unimodal Distributions}
\label{exp:unimodal}

MCMC algorithms which are designed to sample from multimodal distributions are often bespoke for multimodal settings and are less efficient for sampling from unimodal distributions. While the previous experiments have shown that \haram{} is competitive for sampling from multimodal distributions, in this experiment we demonstrate that \haram{} demonstrates improvement over HMC for sampling from unimodal distributions. We consider the problem of sampling from a standard $d$-dimensional Gaussian distribution, $\N(\zerov, \I_d)$, and compare the performance of \haram{} and HMC for $d \in \{2, 10, 50, 100\}$. For both HMC and \haram{}, we set $\epsilon=0.5$ and $L=20$, and we draw samples of size $10,000$ (with the first $5,000$ as warm-up). Additionally, for \haram{} we set $\gamma=0.05$. Therefore, any differences in the performance of the two methods can be purely attributed to the repelling-attracting dynamics~introduced~by~the~parameter~$\gamma$. 

The results in Table~\ref{table:unimodal} present the quality of the last $5,000$ samples obtained from HMC and \haram{} in terms of the $\w$ metric and its expected margin of error  of $\approx n^{-1/d}$ \citep{fournier2015rate} is provided in the last row of Table~\ref{table:unimodal} for reference. Notably, there is no significant difference in sample quality between the two methods. Since $\e$ and $L$ are the same for both methods, the CPU time for both methods is also roughly the same. 

\begin{figure}[t!]
    \includegraphics[width=\textwidth]{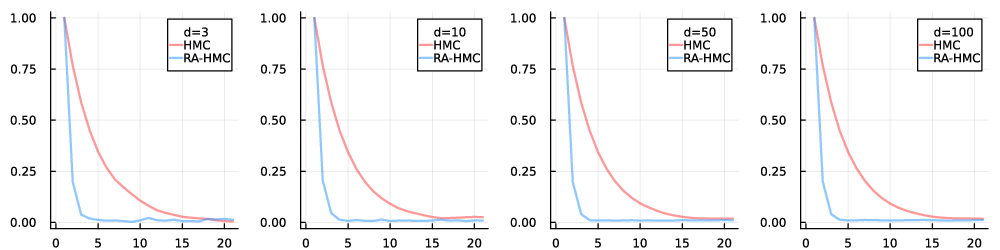}
    \caption{Autocorrelation functions of HMC's and \haram{}'s samples for the standard unimodal, isotropic Gaussian distributions with $d \in
    \{3, 10, 50, 100\}$.}
    \label{fig:unimodal}
\end{figure}

On the other hand, the autocorrelation functions of the last $5,000$ samples (and for the first coordinate) of the Markov chain are shown in Figure~\ref{fig:unimodal}. Notably, in all dimensions \haram{} leads to considerably better mixing than HMC. This can be explained by the fact that, as per \eref{eq:haram-a}, the first half of the \haram{} trajectory is essentially simulating a \textit{chaotic dynamical system} where the modes of the target distribution are Lyapunov unstable states. This leads to better mixing of the samples, and substantially improves the expected sample size (ESS) without sacrificing the sample quality.

\begin{table}[t!]
    \centering
    \caption{$\w$ metric for sampling from unimodal, isotropic Gaussian distribution}\label{table:unimodal}
    \small
    \begin{tabular}{l|rrrrr}\toprule
        Method &$d=3$ &$d=10$ &$d=50$ &$d=100$ \\\midrule
        HMC &0.17 &3.10 &53.41 &133.09 \\
        \haram{} &0.19 &3.23 &54.72 &136.05 \\ \midrule
        Margin of error &0.06 &0.43 &0.84 &0.92 \\
        \bottomrule
    \end{tabular}
\end{table}

%% file: inputs/discussion.tex
% !TEX root = %WORKSPACE_FOLDER%/plain.tex
% \clearpage
\section{Concluding Remarks}
\label{sec:conclusion}

In this work, we have proposed the repelling-attracting Hamiltonian Monte Carlo (\haram{}) as an enhancement of Hamiltonian Monte Carlo for sampling from high-dimensional target distributions with complex, multimodal structures. The key idea which distinguishes \haram{} from HMC is the introduction of the friction parameter, $\gamma$, which facilitates the movement \textit{across} iso-energy contours to generate proposals. The \haram{} algorithm comprises of two stages: in the first stage, the high-density regions of the target distribution serve as repelling states, enouraging mode exploration; while in the second stage, the high-density regions serve as attracting states and forces to trajectory to settle near alternative modes. While this is achieved by departing from conservative Hamiltonian systems to the generalized conformal Hamiltonian systems, \haram{} is designed to perserve the target distribution and the underlying symplectic structure---properties which are key for the success of HMC. We have shown that \haram{} also retains the same numerical accuracy as HMC, while being able to generate Markov chains which lead to substantially better mixing. Lastly, we provided a procedure for automatically tuning the hyperparameters of \haram{} using the version of Nesterov dual-averaging proposed by \cite{hoffman14a}. 

In Section~\ref{sec:extension} we discussed how the repelling-attracting idea underlying \haram{} can be extended to some other variants of HMC in literature. However, the benefit of such an extension has not been investigated in this article, and is an interesting direction for future research. Furthermore, the No-U-Turn Sampler \citep[NUTS,][]{hoffman14a} is a popular variant of HMC which has been shown to be more efficient than HMC in many cases. It would be interesting to investigate the possibility of incorporating the repelling-attracting mechanism into NUTS, and the impact of such an extension. We point out, however, that in stark contrast to NUTS, we want to ensure that the trajectory generated by \haram{} makes a U-Turn when it transitions from the repelling to the attracting stage. Therefore, a direct application of the repelling-attracting mechanism to NUTS may not be straightforward. Lastly, we have not been able to guarantee that \haram{} conserves energy, although the energy drift can be bounded as shown in Proposition~\ref{prop:energy}, and is able to maintain high acceptance probability even across distant trajectories, as evidenced by the numerical experiments in Section~\ref{sec:xperiments}. A more detailed theoretical investigation of the energy preservation and the geometric ergodicity of \haram{} is left for future research. 

%% file: inputs/proofs.tex
\section{Proofs}
\label{sec:proofs}

In this section, we provide the proofs of the main results in the paper.

\subsection{Proof of Proposition~\ref{prop:mode-transition}}
\label{proof:prop:mode-transition}

\begingroup

\renewcommand{\st}{\qq{s.t.}}
\renewcommand{\implies}{\Longrightarrow}

We begin by noting that the randomness in the event $\E(b, d, \alpha)$ arises solely from the initial momentum $\po \sim \N(\zerov, \I)$. With this in mind, let $B$ be the hyperplane
defined by
\eq{
    B = \qty{ \yv \in \R^d : \norm{\yv - \bd} = \norm{\yv + \bd} },
}
and let $B^-$ and $B^+$ be the half-spaces which are closer to $-\bd$ and $\bd$, respectively. Let $\E_1, \E_2$ be the events defined by,
$$
\E_1 = \qty\Big{ \po : \half \norm{\po}^2 \ge \inf_{\q \in \Aa}\log\pi(\q) - \sup_{q \in B} \log\pi(\q) } \qq{and} \E_2 = \qty\Big{ \po : \half \norm{\po}^2 \ge d \qty(\frac{b^2 - \alpha^2}{4\sigma^2})}.
$$
The proof proceeds in two steps:
\begin{enumerate}
    \item We first show that $\E(b, d, \alpha) \subseteq \E_1 \subset \E_2$
    \item We use a Chi-squared tail probability bound on $\pr(\E_2)$ to obtain the desired result.
\end{enumerate}

\textbf{Step 1:} For any $\qo \in \Aa$ and $\po \in \E(b, d, \alpha)$, since $\qo \in B^+$ and $\eta(T; \qo, \po) \in B^-$, and because  $\Phi_t: \R^{2d} \rightarrow \R^{2d}$ is a smooth map for each $t$, it follows that there exists $t \in [0, T]$ such that $\eta(t; \zo) \in B$. Therefore, 
\allowdisplaybreaks[1]
\eq{
    \po \in \E(b, d, \alpha) &\implies \exists t \in [0, T] \st \eta(t; \po, \qo)\in \E_1\nn\\
    &\implies \exists t \in [0, T] \st H(\qp[t]) \le H(\qp[o]), \quad \qt \in B \nn\\
    &\implies \exists t \in [0, T] \st \half \norm{\pt}^2 - \log\pi(\qt) \le \half \norm{\po}^2 - \log\pi(\qo), \quad \qt \in B \nn\\
    &\implies \exists t \in [0, T] \st \half \norm{\po}^2 \ge  \log\pi(\qo) - \log\pi(\qt) + \half \norm{\pt}^2, \quad \qt \in B\nn\\
    &\implies \exists t \in [0, T] \st \half \norm{\po}^2 \ge  \inf_{\q \in \Aa}\log\pi(\q) - \sup_{q \in B} \log\pi(\q) + \half \norm{\pt}^2 \nn\displaybreak\\
    \therefore \po \in \E(b, d, \alpha) &\num{i}{\implies} \norm{\po}^2 \ge  \inf_{\q \in \Aa}\log\pi(\q) - \sup_{q \in B} \log\pi(\q) \nn\\
    &{\implies} \po \in \E_1,
}
where the implication (i) follows from the fact that $\norm{\pt}^2 \ge 0$. Therefore, $\E(b, d, \alpha) \subseteq \E_1$. Furthermore, for $\pi(\q)$ given by 
$$
\pi(\q) = \lambda \qty{\exp\qty(-\frac{1}{2\sigma^2} \norm{\q - \bd}^2) + \exp\qty(-\frac{1}{2\sigma^2} \norm{\q + \bd}^2)},
$$
where $\lambda > 0$ is the normalizing constant of $\pi$, using the fact that $\log(x+y) > \log x$ for all $x, y > 0$, it follows that
\eq{\label{eq:inf-bound}
    \inf_{\q \in \Aa}\log\pi(\q) &\ge \log\lambda - \sup_{\q \in \Aa} \frac{\norm{\q - \bd}^2}{2\sigma^2} = \log\lambda - \frac{\alpha^2d}{2\sigma^2}.
}
Similarly, using the fact that $\norm{\q - \bd} = \norm{\q + \bd}$ for all $\q \in B$, it follows that
\eq{\label{eq:sup-bound}
    \sup_{\q \in B} \log\pi(\q) &\le \log2 + \log\lambda - \inf_{\q \in B} \frac{\norm{\q - \bd}^2}{2\sigma^2} = \log2 + \log\lambda - \frac{b^2d}{2\sigma^2}.
}
For $\po \in \E_1$,  plugging in \eref{eq:inf-bound} and \eref{eq:sup-bound} into lower bound on $\half\norm{\po}^2$ we get
\eq{
    \half \norm{\po}^2
    &\ge \inf_{\q \in \Aa}\log\pi(\q) - \sup_{q \in B} \log\pi(\q)\nn\\
    &\ge \log\lambda - \frac{\alpha^2d}{2\sigma^2} - \log2 - \log\lambda + \frac{b^2d}{2\sigma^2}\nn\\
    &= \frac{(b^2-\alpha^2)d}{2\sigma^2} - \log2\nn\\
    &\num{ii}{>} \frac{(b^2-\alpha^2)d}{4\sigma^2},
}
where (ii) follows from the fact that $b \ge \sqrt{\alpha^2 + 2 \sigma^2} > \sqrt{\alpha^2 + (4\log 2/d)\sigma^2}$ whenever $d \ge 2$. Therefore, we have that $\E(b, d, \alpha) \subseteq \E_1 \subset \E_2$.

\textbf{Step 2:} We now use a Chi-squared tail probability bound on $\pr(\E_2)$ to obtain the desired result. Using the fact that $\norm{\po}^2 \sim \chi^2_d$ when $\po \sim \N(\zerov, \I)$, the Chi-squared tail probability bound yields
\eq{
    \pr(\E(b, d, \alpha)) < \pr\pa{\E_2} &= \pr\pa{ \norm{\po}^2 \ge \f{(b^2-\alpha^2)d}{2\sigma^2} }\nn\\
    &\num{iii}{\le} \exp\qty(-\frac{d}{2} \qty{ \f{b^2-\alpha^2}{2\sigma^2} - 1 - \log\qty(\f{b^2-\alpha^2}{2\sigma^2}) }),\nn
}
where (iii) uses \citet[Eq.~3]{ghosh2021exponential}, which holds whenever $b > \sqrt{\alpha^2 + 2\sigma^2/d}$. This gives us the desired result.\QED
\endgroup

%%%%%%%%%%%%%%%%%%%%%%%%%%%%%%%%%%%%%%%%%%%%%%%%%%%%%%%%%%%%%%%%%%%%%%%%%%%%%%
%%%%%%% Proof: Lemma Energy
%%%%%%%%%%%%%%%%%%%%%%%%%%%%%%%%%%%%%%%%%%%%%%%%%%%%%%%%%%%%%%%%%%%%%%%%%%%%%%

\subsection{Proof of Lemma~\ref{lemma:energy}}
\label{proof:lemma:energy}

\begingroup
For $\z_t = \up_t(\z_0)$ we have that
\eq{
    \dd{t}H(\z_t) 
    &= \Big\langle{\D_\z H(\z_t), \dd{t}\z_t}\Big\rangle\nn\\ 
    &= \D_\z H(\z_t)\tr \qty\Big(\Om \D_\z H(\z_t) + \G \z_t)\nn\\
    &= \D_\z H(\z_t)\tr \Om \D_z H(\z_t) + \D_\z H(\z_t)\tr \G \z_t\nn\\
    &\num{i}{=} \D_\z H(\z_t)\tr \G \z_t\nn\\ 
    &= \gamma \p_t \tr \S\inv \p_t,\nn
}
where (i) follows from the fact that $\D_\z H(\z_t)\tr \Om \D_z H(\z_t) = \zerov$, and the final equality follows by noting that $\D_\z H(\z_t) = (\grad U(\q_t), \S\inv \p_t)\tr$. A similar argument for $\z_t = \dn_t(\z_0)$ gives the final result. \QED

\endgroup

%%%%%%%%%%%%%%%%%%%%%%%%%%%%%%%%%%%%%%%%%%%%%%%%%%%%%%%%%%%%%%%%%%%%%%%%%%%%%%
%%%%%%% Proof: Reversibility
%%%%%%%%%%%%%%%%%%%%%%%%%%%%%%%%%%%%%%%%%%%%%%%%%%%%%%%%%%%%%%%%%%%%%%%%%%%%%%

\subsection{Proof of Proposition~\ref{prop:reversibility}}
\label{proof:prop:reversibility}

\begingroup
\forcecommand{\etaup}{\eta^{+}}
\forcecommand{\etadn}{\eta^{-}}
\forcecommand{\xiup}{\xi^{+}}
\forcecommand{\xidn}{\xi^{-}}
\forcecommand{\zi}{\zeta}

Let $H$ be the Hamiltonian satifying the condition that $H(\q, \p) = H(\q, -\p)$.  Consider the dynamics of $\dn_t$ given by \eref{eq:haram-b}, i.e., 
\eq{
\label{eq:rev-fwd}
{\dd{t}\q_t = \D_\p H(\q_t,\p_t)}%_{\circled{1}}, 
\qq{and} 
{\dd{t}\p_t = -\D_\q H(\q_t, \p_t) - \gamma \p_t}%_{\circled{2}}.
}
Let $\tz = \flip(\z)$ be given by $(\tq,\tp) = (\q, -\p)$, and take $\tt = -t$. Then,
\eq{
\dd{\tt}\tq_t = -\dd{t}\q_t &= -\D_\p H(\q_t, \p_t) = -\D_{\p} H(\tq_t, \tp_t) = \D_{\tp} H(\tq_t, \tp_t),\nn
}
where the final equality follows by noting that $\pt = -\pt$. Similarly,
\eq{
\dd{\tt}\tp_t = \dd{t} \p_t = -\D_\q H(\q_t, \p_t) - \gamma \p_t  = -\D_\tq H(\tq_t, \tp_t) + \gamma \tp_t,\nn
}
which corresponds to the dyanmics of $\up_t$ in \eref{eq:haram-b}. Therefore, for $(\q_t, \p_t) = \dn_t(\q_0, \p_0)$ we have
$$
\up_t\qty\big(\flip(\q_t, \p_t)) = \up_t(\tq_t, \tp_t) = (\q_0, -\p_0) = \flip \circ \dn_{-t}(\q_t, \p_t),
$$
or, equivalently, $\up_t \circ \flip = \flip \circ \dn_\mt$. An identical argument also yields $\dn_t \circ \flip = \flip \circ \dn_\mt$, and the first claim follows by noting that $\flip \circ \flip = \id$. \bigskip

For the second claim, taking $\F_T = \dn_{T/2} \circ \up_{T/2}$, it follows that
\eq{
\qty(\flip \circ \F_T) \circ \qty(\flip \circ \F_T) &=
\qty(\flip \circ \dn_{T/2} \circ \up_{T/2}) \circ \qty(\flip \circ \dn_{T/2} \circ \up_{T/2})\n
&= \flip \circ \dn_{T/2} \circ \up_{T/2} \circ \underbrace{\flip \circ \dn_{T/2}}_{=\up_{-T/2} \circ \flip} \circ \up_{T/2}\n
&= \flip \circ \dn_{T/2} \circ \underbrace{\up_{T/2} \circ \up_{-T/2}}_{=\id} \circ \flip \circ \up_{T/2}\n
&= \underbrace{\flip \circ \dn_{T/2} \circ \flip}_{=\up_{-T/2}} \circ \up_{T/2}\n
&= \up_{-T/2} \circ \up_{T/2} = \id,\nn
}~and, therefore, $\flip \circ \F_T$ is an involution.

For the discrete dynamics, consider one step of the conformal leapfrog step in Algorithm~\ref{alg:haram}, i.e., 
\eq{
\P^{\pm}_\e(\qp) = \qty\Big(\eta^{\pm} \circ \underbrace{\xi_1 \circ \xi_2 \circ \xi_1}_{=:\xi} \circ \eta^{\pm})(\qp),
\label{eq:discrete-dynamics}
}
where 
\eq{
    \eta^{\pm}(\q, \p) := (\q, e^{\pm\gamma\ebyt}\p),\qq{} \xi_1(\q, \p) := (\q, \p - \ebyt \grad U(\q)),\qq{and} \xi_2(\q, \p) := (\q + \e\S^{-1}\p, \p).\nn
}
For $(\tq,\tp) = (\q, -\p)$, observe that
\eq{
    (\etaup \circ \flip \circ \etadn)(\tq, \tp) &= \etaup(\q, -e^{-\gamma\ebyt}\p) = (\q, -\p),\nn
}
which implies that $\etaup \circ \flip \circ \etadn = \flip$. Furthermore, the by noting that the usual leapfrog step $\flip \circ \xi$ is an involution, this implies that
\eq{
    \dn_\e \circ \flip \circ \up_\e 
    &= \etadn \circ \xi \circ \underbrace{\etadn \circ \flip \circ \etaup}_{=\flip} \circ \xi \circ \etaup\nn\\
    &= \etadn \circ \underbrace{\xi \circ \flip \circ \xi}_{=\flip} \circ \etaup\nn\\
    &= \etadn \circ \flip \circ \etaup\nn\\
    &= \flip.\nn
}
where, the second equality follows from the fact that $\flip \circ \xi = \xi \circ \flip$ for the usual leapfrog step. A repeated application of the result above yields
\eq{
    \up_{\el} \circ \flip \circ \dn_{\el} 
    &= \qty(\up_{\e})^L \circ \flip \circ \qty(\dn_{\e})^L = \flip,\nn
}
and, therefore, $\flip \circ \F_\el \circ \flip \circ \F_\el = \id$, which is the desired result.\QED

% \eq{
% }
\endgroup

%%%%%%%%%%%%%%%%%%%%%%%%%%%%%%%%%%%%%%%%%%%%%%%%%%%%%%%%%%%%%%%%%%%%%%%%%%%%%%
%%%%%%% Proof: Symplecticity
%%%%%%%%%%%%%%%%%%%%%%%%%%%%%%%%%%%%%%%%%%%%%%%%%%%%%%%%%%%%%%%%%%%%%%%%%%%%%%

\subsection{Proof of Proposition~\ref{prop:symplecticity}}
\label{proof:prop:symplecticity}

The preservation of symplectic structure and, therefore, volume for the continuous time dynamics is a simple consequence of \citet[Proposition~1]{mclachlan2001conformal}. See also \cite{francca2020conformal}. In particular, for the symplectic $2$-form
\eq{
    \omega_0 = \d\q_0 \wedge \d\p_0 = \sum_{i=1}^d \d q_{0,i} \wedge \d p_{0,i}\,\,,
}
the symplectic $2$-form at time $t$ along the flow $\varphi_t$ is given by the pullback $\omega_t := (\varphi_t)^*\omega_0 = \d\qt \wedge \d\pt.$ From \citet[Proposition~1]{mclachlan2001conformal}, it follows that
\eq{
    (\up_t)^*\omega_0 = e^{+\gamma t}\omega_0, \qq{and} (\dn_t)^*\omega_0 = e^{-\gamma t}\omega_0.\nn
}
Therefore, for the \haram{} dynamics we have
\eq{
    \omega_t = \begin{cases}
        e^{+\gamma t}\omega_0 & \text{ for } 0 \le t \le T/2\\
        e^{-\gamma (t-T/2)}\omega_{T/2} = e^{\gamma \qty(T - t)}\omega_0 & \text{ for } T/2 \le t \le T,
    \end{cases}
}
and, in particular, at time $t=T$ we have $\omega_T = \omega_0$. Furthermore, since the volume form in the extended phase-space $\R^{2d}$ is the exterior product of the symplectic $2$-form, i.e., $\textup{vol} := \omega^{\wedge d}$, from \citet[Chapter~8.38~B, pg.~206]{arnol2013mathematical}, it follows that
\eq{
    \textup{vol}_t = \begin{cases}
        e^{+d\gamma t}\textup{vol}_0 & \text{ for } 0 \le t \le T/2\\
        e^{-d\gamma (T-t)}\textup{vol}_0 & \text{ for } T/2 \le t \le T,
    \end{cases}
}
and $\textup{vol}_T = \textup{vol}_0$, implying that volume is also conserved in the phase-space at time $t=T$ for the \haram{} dynamics.

For the discrete dynamics in \eref{eq:discrete-dynamics}, the determinant of Jacobian for the leapfrog step is $\det({\jac{\xi}}) = 1$. It is easy to verify that $\det(\jac{\eta^{\pm}}) = \exp(\pm\gamma\ebyt)$. By the Jacobian chain rule, 
$$
\det\qty\big(\jac{\P^{\pm}_\e}) = \det\qty\big(\jac{\big( \eta^\pm \circ \xi \circ \eta^\pm \big)}) = \det(\jac{\eta^{\pm}}) \det(\jac{\xi}) \det(\jac{\eta^{\pm}}) = \exp(\pm\gamma\e).
$$
Therefore, for $\F_\el$ it follows that 
$$
\det\qty\Big(\jac{\F_\el}) = \det\qty\Big(\jac{\big(\dn_\el \circ \up_\el)}) = e^{-L\gamma\e} e^{+L\gamma\e}  = 1,
$$
which implies that volume is also conserved in the phase-space for the transition $\z \mapsto \F_\el(\z)$.\QED

%%%%%%%%%%%%%%%%%%%%%%%%%%%%%%%%%%%%%%%%%%%%%%%%%%%%%%%%%%%%%%%%%%%%%%%%%%%%%%
%%%%%%% Proof: Approximation Error
%%%%%%%%%%%%%%%%%%%%%%%%%%%%%%%%%%%%%%%%%%%%%%%%%%%%%%%%%%%%%%%%%%%%%%%%%%%%%%

\subsection{Proof of Proposition~\ref{prop:approximation}}
\label{proof:prop:approximation}

\begingroup
\forcecommand{\et}{\frac{\e}{2}}
\forcecommand{\s}{^*}
\forcecommand{\y}{{\boldsymbol{y}}}
\forcecommand{\elll}{\epsilon,\ell}
\forcecommand{\tel}{{\epsilon,L/2}}
\forcecommand{\Fs}{\F^*}
\forcecommand{\dns}{\P^{-, *}}
\forcecommand{\ups}{\P^{+, *}}

For simplicity let $L = \floor{T/\e}$, and let $X_{\up}$ be the Lie derivative operator associated with the dynamics of \eref{eq:haram-a}, $\up_t$, given by
\eq{
    X_{\up} = \S\inv \p \ddd{\q} - \qty\Big(\grad U(\q) - \gamma \p) \ddd{\p}.\nn
}
We refer the reader to \citep[Chapter~III,~Section~5]{hairer2006geometric} for a detailed treatment of Lie derivatives and their properties. The Lie derivative operator $X_{\up}$ is a first-order differential operator, and, in particular, from \citet[Chapter~III,~Eq.~3.5]{hairer2006geometric} $X_{\up}$ is the generator of the \textit{exact flow} of $\up$ via the identity
$$
\up_t(\z) = \exp(tX_{\up})\id(\z).
$$
Consider the symmetric (Strang) splitting of $X_{\up}$ given by
\eq{
    X_{\up} = \half A + \half B + C + \half B + \half A.\nn
}
where
\eq{
    &A := \gamma \p \ddd{\p}, \qq{}B := -\grad U(\q) \ddd{\p}, \qq{and} C := \S\inv \p \ddd{\q}.\quad \nn
}
For $\e \approx \d{t}$, one step of the conformal leapfrog algorithm for $\up_\e$ as given in \eref{eq:haram-strang} is equivalently written as
\eq{
    \up_\e = \underbrace{\exp(\et A)}_{=\xi_{\e/2}^{\G}} \circ \underbrace{\exp(\et B) \circ \exp(\e C) \circ \exp(\et B)}_{=\P^H_{\e}} \circ \underbrace{\exp(\et A)}_{=\xi^{\G}_{\e/2}} \circ \id.\nn
}
Using the Baker-Campbell-Hausdorff formula, it follows that
\eq{
    \up_\e &\num{i}{=} \exp\qty(\et A) \circ \exp\qty(\e (B + C) + O(\e^3)) \circ \exp\qty(\et A) \circ \id\nn\\
    &\num{ii}{=} \exp\qty\big(\e (A + B + C) + O(\e^3)) \circ \id\nn\\
    &= \exp\qty(\e X_{\up} + O(\e^3)) \circ \id\nn
}
where (i) and (ii) follow from \citep[Chapter~IX,~Eq.~4.4]{hairer2006geometric}, and the final equality follows from the fact that $X_{\up} = A + B + C$. Mutatis mutandis, an analogous argument for the flow $\dn$ yields
\eq{
    X_{\dn} := \S\inv \p \ddd{\q} - \qty\Big(\grad U(\q) + \gamma \p) \ddd{\p},\qq{and} \dn_\e = \exp\qty(\e X_{\dn} + O(\e^3)) \circ \id.\nn
}
Therefore, for the transition $\z \mapsto \F_\el(\z)$, it follows that
\eq{
    \F_\el &= \dn_\tel \circ \up_\tel\nn\\ 
    &= (\dn_\e)^{L/2} \circ (\up_\e)^{L/2}\nn\\
    &\num{iii}{=} \exp\qty(\f{L\e}{2} X_{\up} + O(L\e^3)) \circ \exp\qty(\f{L\e}{2} X_{\dn} + O(L\e^3)) \circ \id\nn\\
    &\num{iv}{=} \exp\qty(\f T2 X_{\up} + O(\e^2)) \circ \exp\qty(\f T2 X_{\dn} + O(\e^2)) \circ \id\nn\\
    &\num{v}{=} \exp\qty(\f T2 X_{\up} + O(\e^2)) \circ \exp\qty(\f T2 X_{\dn}) \circ \id + O(\e^2)\nn\\
    &\num{vi}{=} \exp\qty(\f T2 X_{\up}) \circ \exp\qty(\f T2 X_{\dn}) \circ \id + O(\e^2)\nn\\
    &= \dn_{T/2} \circ \up_{T/2} + O(\e^2)\nn\\ 
    &= \F_T + O(\e^2)\nn
}
where (iii) follows from Gr\"obner's lemma \citep[Chapter~III,~Lemma~5.1]{hairer2006geometric} and the fact that $\exp(M)^L = \exp(LM)$, (iv) uses the substitution $T = L\epsilon$, and (v, vi) follow from the fact that $\exp(\e^2 M) \circ f = \qty(\id + O(\e^2M)) \circ f = f + O(\e^2)$. Therefore, the error in the approximation of $\F_T$ by $\F_\el$ is
\eq{
    \norm{\F_\el(\z) - \F_T(\z)} = O(\e^2), \qq{for} L = O(1/\e) \text{ steps}.\nn
}
Since $U(\q)$ is $\ell$-Lipschitz, it follows that $\abs{H(\z) - H(\z')} \le C \norm{\z-\z'}$ where $C = \max\qty{\ell, \norm{\S\inv}}$, and, we have that
\eq{
    \abs\big{H\qty(\F_\el(\z)) - H\qty(\F_\el(\z))} \le O(\e^2).\nn
}
\QED

\endgroup

%%%%%%%%%%%%%%%%%%%%%%%%%%%%%%%%%%%%%%%%%%%%%%%%%%%%%%%%%%%%%%%%%%%%%%%%%%%%%%
%%%%%%% Proof: Energy Drift
%%%%%%%%%%%%%%%%%%%%%%%%%%%%%%%%%%%%%%%%%%%%%%%%%%%%%%%%%%%%%%%%%%%%%%%%%%%%%%

\subsection{Proof of Proposition~\ref{prop:energy}}
\label{proof:prop:energy}

\begingroup
% \forcecommand{\pt}{\widetilde{\boldsymbol{p}}}
\forcecommand{\pt}{\widehat{\p}}
\forcecommand{\y}{\boldsymbol{y}}
\forcecommand{\ss}{\calo{S}}
\forcecommand{\by}{\overline{\boldsymbol{y}}}
\forcecommand{\yb}{\overline{\boldsymbol{y}}}
\forcecommand{\bz}{\overline{\boldsymbol{z}}}
\forcecommand{\zb}{\overline{\boldsymbol{z}}}
\forcecommand{\bq}{\overline{\boldsymbol{q}}}
\forcecommand{\qb}{\overline{\boldsymbol{q}}}
\forcecommand{\bp}{\overline{\boldsymbol{p}}}
\forcecommand{\pb}{\overline{\boldsymbol{p}}}
\forcecommand{\hh}{\widehat{H}}
\forcecommand{\W}{\mathcal{W}}
\forcecommand{\halfinv}{^{-1/2}}
\newcommand{\s}{\sigma}

\renewcommand{\dot}[1]{%
  \mathchoice%
  {\accentset{\mbox{\large\bfseries .}}{#1}}% Dot for display style
  {\accentset{\mbox{\large\bfseries .}}{#1}}% Dot for text style
  {\accentset{\mbox{\scriptsize\bfseries .}}{#1}}% Dot for script style
  {\accentset{\mbox{\tiny\bfseries .}}{#1}}% Dot for scriptscript style
}

To ease the notation and to make the presentation of the proof consistent with the literature in dynamical systems and control theory (e.g., \citealp{sastry2013nonlinear,khalil2002nonlinear}), we will make use of several notational conventions. For $\z_t, \q_t, \p_t$, wherver the dependence on time is clear from context, we will drop the subscript $t$ and express them as $\z, \p, \q$, and the  (autonomous) conformal Hamiltonian system from \eref{eq:haram-b} by
\eq{
    \begin{pmatrix}
        \dot \q \\
        \dot \p
    \end{pmatrix}
    = 
    \begin{pmatrix}
        \S\inv\p \\
        -\grad U(\q) - \gamma \p
    \end{pmatrix} =: f(\q, \p), \label{eq:autonomous}
}
where $\gamma > 0$, $(\dot \q, \dot \p) := \qty(\dd{t}\q_t,  \dd{t}\p_t)$, and $f(\qp): \R^{2d} \to \R^{2d}$ is the vector field describing the dynamics of the conformal Hamiltonian system. For a function $G: \R^{2d} \to \R$, the time derivative of $G$ along the trajectories $(\qp[t]) = \dn_t(\qp)$ of \eref{eq:autonomous} is given by $\dot G(t) := \grad G(t)\tr f(\q, \p)$.

\begin{figure}[t!]
    \centering
    \includegraphics[width=0.7\textwidth]{\Root/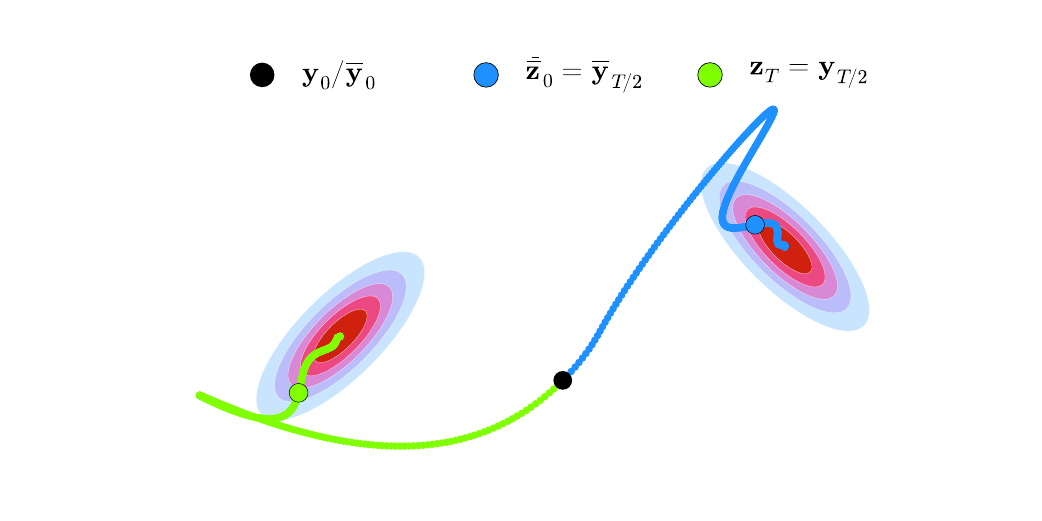}
    \caption{The figure illustrates the key ingredients of the proof in the state-space for $\q \in \R^d$. For $\z_0 \in \R^{2d}$ (green point), we take ${\y_0 := \up_{T/2}(\z_0)}$ and $\yb_0:= \flip(\y_0)$. Since $\flip(\q, \p) = (\q, -\p)$, the projections of $\y_0$ and $\yb_0$ on the state-space $\R^d$ are the same point, shown in black. From Proposition~\ref{prop:reversibility}, it follows that $\y_{T/2} = \dn_{T/2}(\y_0) = \z_T$ and $\yb_{T/2} = \dn_{T/2}(\yb_0) = \z_{0}$ (blue point). By the invariance of the Hamiltonian under $\flip$, $H(\zb_{0}) = H(\z_0)$. We show that $H(\y_t)$ and $H(\yb_t)$ decay exponentially along the green and blue trajectories, respectively. The proof, then, follows by combining these results. We also note that although the black point is the same for $\y_0$ and $\yb_0$ in the state-space; in the phase-space, the momenta are flipped, therefore $\y_0$ and $\yb_0$ are reflections about the $\qty{(\q, \zerov): \q \in \R^d}$ subspace.}
    \label{fig:flip}
\end{figure}

Let $\z_0 \in \R^{2d}$ be the initial state at $t=0$ and let $\z_T = \Psi_T(\z_0)$ be the final proposed state by \haram{} at time $t=T$. We are going to translate the time parameter $t$ to consider $\y_0 = \up_{T/2}(\z_0)$ as the intermediate state of the \haram{} trajectory at time $t=T/2$ and let $\yb_0 = \flip(\y_0)$ be the state obtained by flipping the momentum of $\y_0$. We also $(\uv, \vv) := \yv$ to denote the position and momentum variables for the state $\yv$, so as to distinguish them from $(\q, \p)$ which are used for $\z$. From Proposition~\ref{prop:reversibility}, it follows that $\z_T = \dn_{T/2}(\y_0)$ and $\zb_0 = \dn_{T/2}(\yb_0)$. 

With this background, the outline of the proof is as follows. See Figure~\ref{fig:flip} for an illustration.

\begin{enumerate}
    \item\label{step-1} From assumptions on $U$, we establish that the trajectories $\dn_t(\y_0)$ and $\dn_t(\yb_0)$ are well-defined and unique for all time $t \ge 0$, which pass through $\zb_0$ and $\z_T$. Their limit points\footnote{Technically, their $\omega$-limit points \citep[Definition~5.18]{sastry2013nonlinear}}, $\ell(\y_0) = \lim_{t \to \infty}\dn_t(\y_0)$ and $\ell(\yb_0) = \lim_{t \to \infty}\dn_t(\yb_0)$, are unique points in $\ss$, which are non-degenerate and isolated critical points of $f(\q, \p)$.
    
    \item\label{step-2} Using standard Lyapunov analysis, we show that the Hamiltonian $H(t)$ is exponentially decaying along the trajectory $\dn_t(\y_0)$, i.e.,
    \begin{align}
        (H(t) - H^*) \le (H(0) - H^*) e^{-\lambda t}, \label{eq:lyapunov}
    \end{align}
    where $H^* := H(\ell(\y_0))$ and $\lambda > 0$ is a constant. An identical result holds for the trajectory $\dn_t(\yb_0)$.
    \item\label{step-3} Combining steps~\ref{step-1}~and~\ref{step-2}, we obtain the final bound on $\abs{H(\z_0) - H(\z_T)}$.
\end{enumerate}

\textbf{Step~\ref{step-1}.}\quad We begin by noting that the critical points of $f(\q, \p)$ are the points in $\calo S$, i.e., 
$$
\text{crit} f := \qty{(\q, \p) \in \R^{2d} : f(\q, \p) = \zerov_{2d}} = \qty{(\q^*, \zerov) \in \R^{2d} : \grad U(\q^*) = \zerov_d}.
$$

Moreover, by the assumption that $U$ is a Morse function, it follows that the Hessian $\grad^2 U(\q^*)$ is non-degenrate for all critical points $\q^*$ \citep[Definition~1.7]{nicolaescu2011invitation}, and from the Morse lemma \citep[Lemma~2.2 and Corollary~2.3]{milnor1963morse}, it also follows that the set of critical points $(\q^*, \zerov) \in \calo S$ are isolated.

From assumption \ref{a-2}, $\grad U$ is Lipschitz gradients, and we have
\eq{
    \norm{f(\qp) - f(\q', \p')} 
    &\num{i}\le \sqrt{(1 + \gamma^2)\norm{\p-\p'}^2 + \norm{\grad U(\q) - \grad U(\q')}^2}\nn\\ 
    &\num{ii}\le \sqrt{{(1+\gamma^2)\norm{\p-\p'}^2 + \beta^2 \norm{\q-\q'}^2}}\nn\\
    &\le L_f \norm{(\q, \p )- (\q',\p')},\nn
}
where (i) follows from the triangle inequality, and (ii) invokes assumption~\ref{a-1}, where the constant $L_f = \max\qty{\beta, 1+\gamma}$. It follows that vector field, $f(\q, \p)$, is Lipschitz continuous. Therefore, by the Picard-Lindel\"of theorem, for any initial condition $\y_0 \in \R^{2d}$, there exists a unique solution $\dn_t(\y_0)$ to \eref{eq:autonomous} for all $t \ge 0$ \citep[Theorem3.22]{sastry2013nonlinear}.

For the points $\y_0, \yb_0 \in \R^{2d}$, let $a := H(\y_0) = H(\yb_0)$, $(\uv_t, \vv_t) := \dn_t(\y_0)$ be the trajectory $\y_t$ along \eref{eq:autonomous}, and $H(t) = H(\uv_t, \vv_t)$. From Lemma~\ref{lemma:energy}, it follows that $\dot H(t) = - \gamma \vv_t\tr \S\inv \vv_t \le 0$, i.e., $H(t)$ is non-increasing along the trajectory $\dn_t(\y_0)$. In particular, for the sublevel set of the Hamiltonian, given by $\hh(a) := \qty{(\q, \p) \in \R^{2d}: H(\q, \p) \le a}$, it follows that $(\uv_t, \vv_t) \in \hh(a)$ for all $t \ge 0$. Furthermore, since $U$ is coercive, i.e., $\lim_{\norm{\uv} \to \infty} U(\uv) = \infty$, it follows that the sublevel set $\hh(a)$ is bounded and compact. Applying LaSalle's invariance principle \citep[Section~5.4]{sastry2013nonlinear}, it follows that $\dn_t(\y_0)$ approaches $\ss$ as $t \rightarrow \infty$. An application of \citet[Theorem~4.1]{attouch2000heavy}  further guarantees that $\lim_{t \rightarrow \infty}\dn_t(\y_0) = \ell(\y_0) \in \ss$. An identical result holds for the trajectory $\dn_t(\yb_0)$.

Lastly, for $\z^* \in \ss$, let $\W(\z^*)$ be the set of initial conditions $\z$ such that $\dn_t(\z)$ converges to $\z^*$, i.e., 
\eq{
    \W(\z^*) := \qty{\z \in \R^{2d} : \lim_{t \to \infty}\dn_t(\z) = \z^*}.\nn
}
The set $\W(\z^*)$ is called the domain of attraction/stable manifold of $\z^*$. From \citet[Proposition~5.44]{sastry2013nonlinear}, it follows that $\W(\z^*)$ is an open, invariant set, and, from the previous discussion, it follows that $\W(\ell(\y_0)), \W(\ell(\yb_0)) \subset \hh(a)$ are both open, invariant sets. Therefore, if $\zb_0 \in \W(\z^*)$, it follows that $\y_0 = \dn_{-t}(\zb_0) \in \W(\z^*)$.

To summarize, for $\y_0, \yb_0 \in R^{2d}$ we have that $\ell(\y_0) = \lim_{t \to \infty}\dn_t(\y_0)$ and $\ell(\yb_0) = \lim_{t \to \infty}\dn_t(\yb_0)$ are non-degenerate, and isolated critical points in $\ss$.

\textbf{Step~\ref{step-2}.}\quad We now show that the Hamiltonian $H(t)$ decays exponentially along the trajectory $\dn_t(\y_0)$. We are interested in obtaining an inequality of the type $\dot H(t) \le -\lambda H(t)$ for some $\lambda > 0$. The issue with working with the Hamiltonian, $H(t)$, directly is that $\dot H(t) \le -\gamma \p_t \S\inv \p_t$, suggesting that $H(t)$ may suddenly stop decreasing when $\p_t = \zerov$. This is, however, not the case since the dynamics of \eref{eq:haram} and \eref{eq:autonomous} nudge the momentum away from these singularities. Therefore, in order to establish the exponential decay of $H(t)$, we will work with a modified Lyapunov function, $V(t)$. This technique is standard in the literature on dynamical systems and control theory, and, more recently in optimization. We refer the reader to \cite{begout2015damped,polyak2018optimisation,maddison2018hamiltonian,muehlebach2021optimization,viorel2020asymptotic}, where similar techniques are used.

In particular, for a small $h > 0$ we employ the Lyapunov function
\eq{\label{eq:Vt}
    V(t) := V(\q, \p) = \underbrace{\half \p\tr\S\inv \p + U(\q)}_{=H(\q, \p)} - U(\q^*) + h \cdot \underbrace{\p\tr \S\inv \grad U(\q)}_{=:\eta(\q, \p)\phantom{\half}},
}
which is similar to the Lyapunov function in \citet[Eq.~6]{polyak2018optimisation}. The time derivative of $V(t)$ along $f(\q, \p)$ is given by
\eq{
    \dot V(t) &= \grad V(t)\tr f(\q, \p) \nn\\
    &= \dot H(t) + h \cdot \qty(\grad\eta(\q_t, \p_t)\tr f(\q, \p)) \nn\\
    &= -\gamma \p\S\inv \p + h \cdot \qty\Big(
        \p \tr \qty\big(\S\inv \grad^2 U(\q) \S\inv) \p - \gamma \cdot \p\tr \S\inv \grad U(\q) - \grad U(\q)\tr \S\inv \grad U(\q)
    )\nn\\
    &= - \pt \tr \qty\big(\gamma \I - h\Av) \pt - h\gamma \pt\tr \S\halfinv \grad U(\q) - h\grad U(\q)\tr \S\inv \grad U(\q),\label{eq:lyapunov-1}
}
where $\pt := \S^{-1/2}\p$ and $\Av := \S\halfinv \cdot \grad^2 U(\q) \cdot \S\halfinv$. 
% For $\lambda < \gamma$ to be specified later, we can write $\dot V(t)$ as
% \eq{
%     \dot V(t) &= - \pt\tr \qty\big(\gamma \I - h\Av) \pt - h(\gamma - \lambda) \pt\tr\S\halfinv\grad U(\q) - h\lambda \p\tr \S\inv \grad U(\q) - h\grad U(\q)\tr \S\inv \grad U(\q).\nn
% }
By Young's inequality, the second term can be bounded by
\eq{
    \pt\tr \S\halfinv\grad U(\q) = \pt\tr\qty(\S\halfinv \grad U(\q)) \ge -\frac{1}{2} \pt\tr\pt - \frac{1}{2} \grad U(\q)\tr \S\inv \grad U(\q).
    \label{eq:youngs}
}
Substituting this back into \eref{eq:lyapunov-1}, we obtain
\eq{
    \dot V(t) &\le -\frac{1}{2} \pt\tr\Bv\pt - {h} \grad U(\q)\tr \S\inv \grad U(\q),\label{eq:lyapunov-2}
}
where $\Bv := \qty(\gamma - \frac{h}{2})\I - h\Av$. Using $\s_{\max}(\Xv) = \s_1(\Xv) \ge \s_2(\Xv) \ge \dots \ge \s_{2d}(\Xv) = \s_{\min}(\Xv)$ to denote the eigenvalues of a matrix $\Xv$, the first term in \eref{eq:lyapunov-2} is bounded by
\eq{
    \pt\tr\Bv\,\pt \ge \s_{\min}(\Bv) \norm{\pt}^2 = \s_{\min}\qty(\qty\Big(\gamma - \frac{h}{2})\I - h\Av) \norm{\pt}^2 \ge \qty(\gamma - \frac{h}{2} - h \frac{\varpi}{\s_{\min}(\S)}) \norm{\pt}^2,\nn
}
where the last inequality follows from the fact that $\s_{\max}\pa{\Av} \le \s_{\max}\pa{\S\inv}\s_{\max}\pa{\grad^2 U(\q)} \le \varpi/\s_{\min}(\S)$. Therefore, $\pt\tr\Bv\,\pt \ge 0$ whenever $h \in (0, h_0)$ for
$$
h_0 := \frac{2\gamma \s_{\min}(\S)}{2\varpi + \s_{\min}(\S)}.
$$
For simplicity, we will take $h = h_0/2$, which implies that $\pt\tr\Bv\pt \ge \gamma \norm{\pt}^2/2$. Similarly, the second term in \eref{eq:lyapunov-2} is bounded by
\eq{
    \grad U(\q)\tr \S\inv\grad U(\q) 
    \ge \frac{1}{\s_{\max}(\S)} \norm{\grad U(\q)}^2 \ge \frac{\mu}{\s_{\max}(\S)} \pa{U(\q) - U(\q^*)},\nn
}
where the last inequality follows from the fact that $U$ satisfies the \pl{} condition. Plugging this into \eref{eq:lyapunov-2}, we obtain
\eq{
    \dot V(t) \le - \frac{\gamma}{2}\norm{\pt}^2 - \frac{h_0\mu}{2\s_{\max}(\S)} \qty\Big(U(\q) - U(\q^*)) \le - \lambda_1 V(t),\label{eq:lyapunov-3}
}
for 
$$
\lambda_1 := \min\qty{\frac{\gamma}{2}, \frac{h_0\mu}{2\s_{\max}(\S)}}.
$$
Using the Gr\"onwall-Bellman lemma \citep[Proposition~3.21]{sastry2013nonlinear}, it follows that
\eq{
    V(t) \le V(0) e^{-\lambda_1 t}.\label{eq:lyapunov-4}
}
Similarly, taking $W(t) = U(\q) - U(\q^*)$, it follows that
\eq{
    \dot W(t) &= \grad W(\q, \p) \tr f(\q, \p)\nn\\
    &= \p\tr \S\inv \grad U(\q)\nn\\
    &\num{iii}{=} \frac{1}{h} \cdot \qty({V(t) - \half\norm{\pt}^2 - W(t)})\nn\\
    &\num{iv}{\le} -\frac{1}{h}W(t) + \frac{1}{h} V(0)e^{-\lambda_1 t},
}
where (iii) follows from the definition of $V(t)$ in \eref{eq:Vt}, and (iv) follows from \eref{eq:lyapunov-4} and by noting that $\norm{\pt}^2 \ge 0$. Therefore, by the method of variation of parameters it follows that $W(t)$ decays exponentially as well with a rate $\lambda_2 = \min\qty{\lambda_1, 1/h}$.

Lastly, for $X(t) = \half \p\tr\S\inv\p$, it follows that
\eq{
    \dot X(t) 
    &= \grad X(\q, \p) \tr f(\q, \p)\nn\\
    &= -\gamma \p\tr\S\inv\p - \p\tr\S\inv\grad U(\q)\nn\\
    &\num{v}{\le} -\gamma \pt\tr\pt - \half \pt\tr\pt - \frac{1}{2}{\grad U(\q)}\S\inv\grad U(q)\nn\\
    &\num{vi}{\le} -\qty(\gamma + \half) \norm{\pt}^2\nn\\
    &= -\qty(2\gamma + 1) X(t) \nn
}
where (v) follows from Young's inequality in \eref{eq:youngs} and (vi) follows from the the fact that $\grad U(\q) \S\inv \grad U(\q) \ge 0$. Therefore, $X(t)$ decays exponentially as well with a rate $\lambda_3 = 2\gamma + 1$, i.e., $X(t) = X(0) e^{-\lambda_3 t}$. Combining the results we get
\eq{
    H(t) - H^* &= X(t) + W(t)\nn\\
    &\le X(0) e^{-\lambda_3 t} + W(0) e^{-\lambda_2 t}\nn\\
    &\le (X(0) + W(0)) e^{-\min\qty{\lambda_2, \lambda_3} t}\nn\\
    &= (H(0) - H^*) e^{-\lambda t},\label{eq:lyapunov-final}
}
where 
\eq{
\lambda = \min\qty{\frac{\gamma}{2}, \frac{\gamma\mu\s_{\min}(\S)}{(2\varpi + \s_{\min}(\S))\s_{\max}(\S)}, \frac{2\varpi + \s_{\min}(\S)}{2\gamma \s_{\min}(\S)}}.
\label{eq:lambda}
}

\textbf{Step~\ref{step-3}.}\quad Let $\z_0 \in \R^{2d}$ and $\z_T = \F_T(\z_0)$, and as described earlier let $\y_0 = \up_{T/2}(\z_0)$ and $\yb_0 = \flip(\y_0)$. Furthemore, let $\y^* = \ell(\y_0)$ and $\yb^* = \ell(\yb_0)$ be the critical points that the trajectories starting from $\y_0$ and $\yb_0$ at time $t=0$ converge to along $\dn_t$. From the results in Step~\ref{step-2}, it follows that
\eq{
    \abs{H(\z_0) - H(\z_T)} 
    &\num{vii}{=} \abs\Big{H(\zb_0) - H(\yb^*) + H(\yb^*) - H(\y^*) + H(\y^*) - H(\z_T)}\nn\\
    &\num{viii}{=} \abs\Big{H(\yb_{T/2}) - H(\yb^*)} + \abs\Big{H(\yb^*) - H(\y^*)} + \abs\Big{H(\y^*) - H(\y_{T/2})}\nn\displaybreak\\
    \therefore \abs{H(\z_0) - H(\z_T)}  &\num{ix}{\le} \qty(H(\yb_0) - H(\yb^*)) e^{-\lambda T/2} + \qty(H(\y_0) - H(\y^*)) e^{-\lambda T/2} + \abs\Big{H(\yb^*) - H(\y^*)}\nn\\
    &\num{x}{=} \qty\Big(H(\z_{T/2}) - H(\y^*) - H(\yb^*)) e^{-\lambda T/2} + \abs\Big{H(\yb^*) - H(\y^*)},\nn
}
where (vii) uses the fact that $H(\z_0) = H(\zb_0)$, (viii) uses the fact that $\z_T = \dn_{T/2}(\y_0)$ and $\zb_0 = \dn_{T/2}(\yb_0)$, (ix) uses \eref{eq:lyapunov-final}, and (x) uses the fact that $\y_0 = \up_{T/2}(\z_0)$ and $H(\y_0) = H(\yb_0)$. Taking the supremum over all $\y^*, \yb^* \in \ss$,
\eq{
    \abs{H(\z_0) - H(\z_T)} \le \sup_{\av, \bv \in \ss}\qty\Big(H(\z_{T/2}) - H(\av) - H(\bv)) e^{-\lambda T/2} + \abs\Big{H(\av) - H(\bv)},\nn
}
which gives us the final result.\QED

\endgroup